\documentclass[final,reqno,oneeqnum,onefignum,onetabnum]{siamart190516}

\usepackage{amssymb}
\usepackage{mathtools}
\usepackage{booktabs}
\usepackage{hyperref}
\usepackage{multirow}
\usepackage{microtype}

\usepackage[T1]{fontenc}
\usepackage{lmodern}

\title
{Self-accelerating root search and optimisation methods based on rational interpolation}

\ifpdf
\hypersetup{ % Optional PDF information
	pdftitle={Self-accelerating root search and optimisation methods based on rational interpolation},
	pdfauthor={S. Cassel}
}
\fi

\headers
{Self-accelerating root search and optimisation methods}
{S.~Cassel}

\author{Sebastian~Cassel\thanks{BNP Paribas, 10 Harewood Avenue, London, NW1 6AA, UK 
		(\email{sebastian.cassel@protonmail.com}).}}

\begin{document}
\maketitle
\thispagestyle{empty} % removes page number

\begin{abstract} %115 words
Iteration methods based on barycentric rational interpolation are derived that exhibit accelerating orders of convergence. %15 words
For univariate root search, the derivative-free methods approach quadratic convergence and the first-derivative methods approach cubic convergence. %19 words (compound=2)
For univariate optimisation, the convergence order of the derivative-free methods approaches $1.62$ and the order of the first-derivative methods approaches $2.42$. %23 words (number=1)
Generally, performance advantages are found with respect to low-memory iteration methods. %12 words
In optimisation problems where the objective function and gradient is calculated at each step, the full-memory iteration methods converge asymptotically $1.8$~times faster than the secant method. %29 words
Frameworks for multivariate root search and optimisation are also proposed, though without discovery of practical parameter choices. %17 words
\end{abstract}


% automatic preconditioning, optimisation, scattered data interpolation

\begin{AMS}
%https://mathscinet.ams.org/mathscinet/msc/msc2010.html
26C15,	%Real rational functions
41A20,	%Approximation by rational functions
65B99,	%Acceleration of convergence in numerical analysis
65D25,	%Numerical differentiation
65F08,	 %Preconditioners for iterative methods
65H05,	%Numerical computation of solutions to single equations
65H10	%Numerical computation of solutions to systems of equations
\end{AMS}
%65H04,	%Numerical computation of roots of polynomial equations

%ACM computing classification
%https://www.acm.org/publications/computing-classification-system/1998/g.1
%G.1.0;	%General
%G.1.1;	%Interpolation
%G.1.2;	%Approximation
%G.1.4;	%Quadrature and Numerical Differentiation
%G.1.5;	%Roots of Nonlinear Equations ***
%G.1.6;	%Optimization ***

% basin of attraction
% complex functions
% extend Brent's method
% tune leading errors of specific problems
% transformation for error disappearance
% other classes of interpolation functions
% line integral (or non-perturbative series solution)
% direct calculation of higher order iteration function
% analytically simple derivative, perhaps sparse, perhaps vanishes (polynomial)
% adjoint differentiation or quantum amplitude estimation
% Grassman numbers

\vspace{1mm}
\phantom{.}

%\vfill
\vspace{4mm}
\section{Introduction}

% introduce families/class of iteration methods
% identify suitable choices for free parameters
% convergence properties

The search for solutions of non-linear equations and related problem of optimisation form routine tasks in numerical analysis. In this paper, families of iterative algorithms are presented that re-use information from an arbitrary number of prior steps in order to accelerate convergence. Furthermore, the iteration parameters are analytically expressed and therefore need no intermediate calibration.

On the history of related work, the general principles of interpolation-based iteration methods for root search problems were studied in the 1960s by Ostrowski~\cite{ostrowski1960solution} and Traub~\cite{traub1964iterative}. That same decade, Tornheim~\cite{tornheim1964convergence}, Jarratt and Nudds~\cite{jarratt1965use} proposed the application of rational interpolation for root searches. However, the general formulations required matrix equations to be solved. In 1980, Larkin~\cite{larkin1980root} then identified a recursive approach for evaluating root search iterations derived from minimal rational interpolants with a linear numerator. Following an alternative approach, in 2008 Sidi~\cite{sidi2008generalization} presented a root search scheme applying Newton's method to polynomial interpolants, with the interpolant derivative being calculated by recursion. These recursive schemes can suffer from numerical instabilities though, depending on the distribution and ordering of interpolation points~\cite{higham2004numerical}.

This paper presents iteration methods based on barycentric representations of rational interpolants~\cite{schneider1986some, berrut2004barycentric}, which possess useful properties of expressivity, computational efficiency and numerical stability~\cite{higham2004numerical, webb2012stability}. Applications for both root search and optimisation problems are considered. Frameworks for multivariate cases are also proposed. However, only the univariate methods are of immediate practical use. It is still common though to perform one-dimensional line searches as a subroutine within multivariate root search or optimisation, and so these methods can be applied to many general problems. 

% The methods offer performance advantages with respect to standard low-memory iteration methods.
% The associated class of non-minimal rational interpolants is thus assessed. 

% key results -> simple algorithms, direct calculation, high efficiency, arbitrary order
% why rational -> free parameters, barycentric rep easy, O(n) update, numerically stable
% special case sub-class algorithms previously know -> generalisation

% simple zeros of real (scalar) function, one-point (greedy) methods
% arbitrary memory, accelerating parameters simple analytic formula
% iteration algorithms with memory 

\vfill
\pagebreak
\section{Roots of univariate rational interpolants}

% evaulation in $\mathcal{O}[n]$ operations
% ease of adding (x,f)
% numerical stability
% don't need to solve any interpolation parameters, $\omega_i^{} \neq 0$, $x_i^{}$ distinct

% invertible in local neighbourhood of simple root
% explain why formula trivially interpolates the function
% iteration formula can be applied to an arbitrary selection of points

The simplest approach for solving $f[x]=0$ is to approximate the inverse function~$x[f]$ and then evaluate at the point $f=0$. Given a set of points $\{ x_i^{}, f_i^{} \}$ where $f_i^{} = f[x_i^{}]$, the barycentric representation of a rational interpolant is stated below,
\begin{align}
x[f] ~\approx\,
\left.
\left( \, \sum_{i=0}^{n} \, \frac{\omega_i^{} \, x_i^{}}{f-f_i^{}} \, \right) \! \middle/ \! \left( \, \sum_{i=0}^{n} \, \frac{\omega_i^{}}{f-f_i^{}} \, \right) \right.
\label{eq:x[f]}
\end{align}

\noindent
where $\omega_i^{}$ are free (non-zero) parameters. The iteration scheme below then follows. 
\begin{align}
\boxed{\hspace{1mm}
	\begin{array}{c}
	\\[-3mm]
	\displaystyle
	x_{n+1}^{} \hspace{1mm}=\, 
	\left.
	\left( \, \sum_{i=0}^{n} \, \frac{\omega_i^{} \, x_i^{}}{f_i^{}} \, \right) \! \middle/ \! \left( \, \sum_{i=0}^{n} \, \frac{\omega_i^{}}{f_i^{}} \, \right) \right.
	\\[4mm]
	\end{array}
}\hspace{1mm}\phantom{.}
\label{eq:method0}
\end{align}

\noindent
Now, instead of addressing whether the inverse function is well-defined or not, it's useful to note that the same iteration formula can also be derived by approximating $f[x]$. 
Consider the rational expression in \cref{eq:f[x]}, where $r[x]$ is approximated by \cref{eq:r0}.
\begin{align}
	f[x] ~=~ \frac{x-c}{r[x]}
\label{eq:f[x]}
\end{align}
\begin{align}
	r[x] ~\approx\,
	\left.
	\left( \, \sum_{i=0}^{n} \, \frac{\omega_i^{} \, r_i^{}}{x-x_i^{}} \, \right) \! \middle/ \! \left( \, \sum_{i=0}^{n} \, \frac{\omega_i^{}}{x-x_i^{}} \, \right) \right.
	\label{eq:r0}
\end{align}

\noindent
From \cref{eq:f[x]}, it is known that $r_i^{} = (x_i^{}-c)/f_i^{}$. The root of the approximated $f[x]$ is then given by \cref{eq:method0} when applying the following constraint and solving for $c$.
\begin{align}
	 \sum_{i=0}^{n} \, \omega_i^{} \, r_i^{} \,=\, 0
	 \label{eq:rfix}
\end{align}

\noindent
Other constraints could be considered, or the right hand side of \cref{eq:rfix} set to some free parameter, leading to alternative iteration formulae. However, this paper does not seek to exhaustively discuss all such possible iteration schemes. It is chosen to focus on the methods that form a weighted linear average of points, as in \cref{eq:method0}.

Before exploring how to choose the interpolation parameters, the case where first derivatives $f_i^\prime$ are also available is presented. The corresponding barycentric rational approximation of the inverse function can be found by taking the limit of \cref{eq:x[f]} for coinciding pairs of interpolation nodes,
\begin{align}
x[f] ~\approx\,
\left.
\left( \, \sum_{i=0}^{n} \, \frac{\lambda_i^{} \, x_i^{} + (\gamma_i^{} \, x_i^{} + \lambda_i^{} / f_i^\prime) \, (f-f_i^{}) }{(f-f_i^{})^2} \, \right) \! \middle/ \! \left( \, \sum_{i=0}^{n} \, \frac{\lambda_i^{} + \gamma_i^{} \, (f-f_i^{})}{(f-f_i^{})^2} \, \right) \right.
\label{eq:xfHermite}
\end{align}

\noindent
where $\{\lambda_i^{}, \gamma_i^{} \}$ are free parameters. The iteration scheme below 
% in eq~\ref{eq:method1} 
then follows for $f=0$.
\begin{align}
\boxed{\hspace{1mm}
	\begin{array}{c}
	\\[-3mm]
	\displaystyle
	x_{n+1}^{} ~=~
	\left.
	\left( \, \sum_{i=0}^{n} \, \frac{ \lambda_i^{} \, ( x_i^{}  - f_i^{} / f_i^\prime ) - \gamma_i^{} \, f_i^{} \, x_i^{} }{ f_i^2 } \, \right) \! \middle/ \! \left( \, \sum_{i=0}^{n} \, \frac{ \lambda_i^{}  - \gamma_i^{} \, f_i^{} }{f_i^2} \, \right) \right.
	\\[4mm]
	\end{array}
}
\label{eq:method1}
\end{align}

\noindent
The same iteration formula can again be derived from $f[x]$ in \cref{eq:f[x]}, but now with:
\begin{align}
r[x] ~\approx\,
\left.
\left( \, \sum_{i=0}^{n} \, \frac{\tilde{\lambda}_i^{} \, r_i^{} + (\gamma_i^{} \, r_i^{} + \tilde{\lambda}_i^{} \, r_i^\prime) \, (x-x_i^{}) }{(x-x_i^{})^2} \, \right) \! \middle/ \! \left( \, \sum_{i=0}^{n} \, \frac{\tilde{\lambda}_i^{} + \gamma_i^{} \, (x-x_i^{})}{(x-x_i^{})^2} \, \right) \right.
\label{eq:r1}
\end{align}

\noindent
and the following constraint, 
\begin{align}
\sum_{i=0}^{n} \, (\gamma_i^{} \, r_i^{} + \tilde{\lambda}_i^{} \, r_i^\prime) \,=\, 0
\label{eq:rfix1}
\end{align}

%\begin{eqnarray}
%r_i^\prime &=& \frac{f_i^{} - (x_i^{}-c) f_i^\prime}{f_i^{2}}
%\end{eqnarray}

\noindent
where the relation $r_i^\prime = (f_i^{} - (x_i^{}-c) f_i^\prime) / f_i^{2}$ is known by differentiating \cref{eq:f[x]}. 
On solving for $c$ and defining $\lambda_i^{} = \tilde{\lambda}_i^{}\, f_i^{\prime}$, the iteration formula in \cref{eq:method1} is reproduced. 

The constraints in \cref{eq:rfix} and \cref{eq:rfix1} could otherwise be identified by first expressing the $r[x]$ approximations as a ratio of polynomial functions: multiplying the numerators and denominators by $\prod_{i=0}^{n} \, (x-x_i^{})^m$ with $m=1,2$ respectively. The constraints are then equivalent to setting the coefficient of the highest degree term $x^{m(n+1)-1}$ in the numerator to zero. Again, other constraints may be considered, but for this choice the iteration scheme corresponds to a weighted linear average of the sampled points and respective Newton iterates.

For future reference, it is also helpful to stress that polynomial functions form a sub-class of barycentric rational representations. The relevant interpolation parameters can be determined from the following partial fraction expansions:
\begin{equation}
\begin{array}{lclll}
\displaystyle \prod_{i=0}^{n} ~ \frac{1}{z-z_i^{}} \!\! &=& \!\!
\displaystyle \sum_{i=0}^{n} \, \frac{u_i^{}}{z-z_i^{}} & \hspace{2mm}
\text{where} \hspace{2mm}
\displaystyle u_i^{} ~=~ \prod_{j\neq i}^{} \, \frac{1}{z_i^{}-z_j^{}}
\\[7mm]
\displaystyle \prod_{i=0}^{n} \, \frac{1}{(z-z_i^{})^2} \!\! &=& \!\!
\displaystyle \sum_{i=0}^{n} \, \frac{u_i^{2} + v_i^{} \, (z-z_i^{})}{(z-z_i^{})^2} & \hspace{2mm}
\text{where} \hspace{2mm}
\displaystyle v_i^{} = -2\, u_i^{2} \, \sum_{j\neq i}^{} \, \frac{1}{z_i^{}-z_j^{}}
\end{array}
\label{eq:pfe}
\end{equation}

\noindent
The approximations for $r[x]$ or $x[f]$ thus simplify to polynomial expressions when the respective denominators equate to products of inverse polynomials.

To determine convergence properties of the iterative methods, series expansions can be considered about the true root $x_\ast^{}$ where $x_i^{} = x_\ast^{} + \epsilon_i^{}$. For the applications presented in this paper, the derivative-free formula in \cref{eq:method0} is focused on, as corresponding parameters for \cref{eq:method1} will be formed from the limit of coinciding pairs of sampling points. 

For derivative-free iteration \cref{eq:method0}, the leading error is:
\begin{align}
\epsilon_{n+1}^{} ~=\,
%x_\ast^{} +  
\left.
\left( \, \sum_{i=0}^{n} \, \omega_i^{} \, \right) \! \middle/ \! \left( \, \sum_{i=0}^{n} \, \frac{\omega_i^{}}{\epsilon_i^{} }  \, \right) \right.
%+ \mathcal{O}[\epsilon_i^2]
 +\,  \text{higher-order terms}
\end{align}

\noindent
Given arbitrary choices of $\omega_i^{}$, the leading error is thus $\mathcal{O}[\epsilon]$ and so the iterations are only linearly convergent (if they do converge). However, the above error term vanishes when $\sum_i^{} \, \omega_i^{} = 0$. If then successively requiring the leading error term to vanish, the following constraints are established:
\begin{align}
\sum_{i=0}^{n} \, \omega_i^{} \, \epsilon_i^{k} \,=\, 0
\hspace{15mm} 
%0 \leq m \leq n-1
0 \leqslant k\leqslant n-1
\label{eq:constraints}
\end{align}

\noindent
Such constraints correspond to finding the kernel of the Vandermonde matrix, whose solution~\cite{berrut1997matrices} is given by:
\begin{align}
\omega_i^{} ~\propto~
\prod_{j\neq i}^{} \, \frac{1}{\epsilon_i^{} - \epsilon_j^{}}
\label{eq:limitSolution}
\end{align}

\noindent %in the expression for $\omega_i^{}$
Although the $\epsilon_i^{}$ terms are not known directly, it can be chosen to substitute $(\epsilon_i^{} - \epsilon_j^{})$ consistently by $(x_i^{} - x_j^{})$ or $(f_i^{} - f_j^{})$ without re-introducing the eliminated error terms. From \cref{eq:pfe}, it then follows that if the constraints in \cref{eq:constraints}  are satisfied, the leading error for the derivative-free iteration scheme has the following form with $m=1$:
\begin{align}
\epsilon_{n+1}^{} ~=~
E_{m,n+1}^{} \times \Bigg( \, \prod_{i=0}^{n} \, \epsilon_i^{m} \, \Bigg)
 +\, \text{sub-dominant terms}
\label{eq:error}
\end{align}

\noindent
where $E_{m,n+1}^{}$ is a factor composed of $k^{\text{th}}$ derivatives $f_\ast^{(k)}$ evaluated at the true root. For convenience in later expressing $E_{m,n+1}^{}$, the following notation is introduced:
\begin{align}
f_\ast^{(i,j,\ldots, k)} ~=~
f_\ast^{(i)}\, %\times 
f_\ast^{(j)}\,  %\times 
\cdots\, %\times 
f_\ast^{(k)} 
\label{eq:notationfi}
\end{align}

\noindent
The error relation trivially generalises to cases where $m$ interpolation nodes coincide. 
In order to then identify the asymptotic order of convergence, the conventional definition that $\epsilon_{n+1}^{}  \sim \epsilon_n^\ell$ can be recursively applied to \cref{eq:error} to find $ \ell^{n+1}  = m \sum_{i=0}^{n} \ell^i $. Equivalently, the geometric sum can be re-expressed to equate $\ell = (m+1) -m\, \ell^{-(n+1)}$. \Cref{tbl:rootConvergence} presents numerical solutions of $\ell$ for the cases where $m=1$ and $2$. 
The order of convergence accelerates as $n$ increases, tending to $m+1$. However, a long memory is not needed to raise the convergence order close to its asymptotic limit. In practice, one may then be satisfied with methods that use a limited history of points.
\begin{table}[!htbp]
	\centering
	%\tiny
	%\scriptsize
	%\footnotesize
	\small
	%\normalsize
	\begin{tabular}{p{18mm}p{14mm}ccccc}
		\\[1mm]
		\toprule
		%Interpolation data & Error & 1-point & 2-point & 3-point & 4-point & 5-point \\ \midrule
		Data & Error & $n=0$ & $n=1$ & $n=2$ & $n=3$& $n=4$\\ \midrule
		$\{ x_i^{} , f_i^{} \}$ & $\prod_{i=0}^{n} \epsilon_i^{}$ &
		1.00000 & 1.61803 & 1.83929 & 1.92756 & 1.96595 \\[1mm]
		$\{ x_i^{} , f_i^{}, f_i^\prime \}$ & $\prod_{i=0}^{n} \epsilon_i^{2}$ & 
		2.00000 & 2.73205 & 2.91964 & 2.97445 & 2.99165 \\ \bottomrule
	\end{tabular}
	\vspace{2mm}
	\caption{Convergence indices ($\ell : \epsilon_{i+1}^{} \sim \epsilon_i^\ell$) for iteration methods with a leading error term given by \cref{eq:error}.}
	\label{tbl:rootConvergence}
\end{table} % to five decimal places

It is stressed that the relative efficiency of derivative-free methods and higher-derivative methods depends on the computational costs of obtaining higher derivatives. If each derivative takes a similar or longer time as $f[x]$ to be calculated, the derivative-free methods are generally most efficient. Even for low-memory iterations, note that the order of convergence is~$2$ for one Newton step given $(f_i^{}, f_i^\prime)$, but~$2.62$~($1.62^2$) for two secant steps given $(f_i^{}, f_{i+1}^{})$ and the respective memory. As $n \to \infty$, the order of convergence is~$3$ for one first-derivative method step, but~$4$~($2^2$) for two derivative-free method steps. 

Conversely, if $f[x]$ is defined by an integral, or $f[x]$ is known to satisfy a differential equation (which could itself be differentiated), or the derivatives assume simple forms, it may be relatively quick to calculate derivatives. In such cases, interpolation-based schemes that minimise the number of references to $f[x]$ could be favourable. The root search schemes detailed in this paper assume that $f[x]$ is determined at every iteration point, but the interpolation-based approach can be adapted to generally consider mixed interpolation conditions.

\vfill\pagebreak
\Cref{tbl:barycentric0,tbl:barycentric1} present root search schemes with specific parameter choices, and the associated leading error factors $E_{m,n+1}^{}$ for some cases of $n$. The 2-point ($n=1$) derivative-free methods in \cref{tbl:barycentric0} are equivalent to the secant method. The 1-point ($n=0$) first-derivative methods in \cref{tbl:barycentric1} are equivalent to Newton's method. The scheme in the left column of \cref{tbl:barycentric0} is also recognised to result in the same iteration steps formulated by Larkin~\cite{larkin1980root}, though without recursive evaluations of divided differences being required. Some of the leading error factors could be cancelled by forming weighted averages of the different iteration schemes. However, the coefficients of the highest degree factors are common for given $n$, and so such averaging does not generally raise the order of convergence further. It is important to note that the $x_i^{}$ or $f_i^{}$ values referenced in the weights must be distinct. Otherwise some weights become infinite.

%weights can be arbitrarily rescaled without affecting the iteration step. 
%$f_n^{}$ latest (/smallest) term

%The leading error factors for the $x_i^{}$-dependent weights involve smaller coefficients, and 
%$E_{1,2(n+1)}^{} = E_{2,n+1}^{}$

\Cref{tbl:example0,tbl:example1} demonstrate the convergence behaviours of different iteration sequences. %where the progression in the convergence phase is related to the respective order of convergence. 
Although the iterations are presented to high precision, there are few practical cases where such precision is needed. If tolerant to 64-bit machine errors ($\sim 10^{-16}$), the interpolation-based methods save a few steps compared to the secant method or Newton's method. Such savings can still form a useful performance advantage though. The application of non-local information can also further help to enter the convergence phase faster, where robust derivative estimates smooth over localised oscillations.

Returning to the topic of parameter constraints in derivative-free schemes, it was previously argued for the terms of $\mathcal{O}[\epsilon_i^k]$ to cancel each other for $k < n$. However, it was then identified that the derivative-free schemes exhibit sub-quadratic convergence. Any $\mathcal{O}[\epsilon_n^k]$ error terms for $k \geqslant 2$ are thus sub-dominant. The constraints in \cref{eq:generalConstraints} can then be generally applied without compromising the order of convergence, where $g_{j}^{} = \mathcal{O}[\epsilon_n^{}]$ and $h_{j}^{} = \mathcal{O}[\epsilon_n^{}]$.
%[\epsilon_n^{}]
%
\begin{align}
\begin{pmatrix}
1 & \cdots & 1 & 1 \\[1mm]
\epsilon_0^{} & \cdots & \epsilon_{n-1}^{} & g_1^{} \\[1mm]
\vdots & & \vdots & \vdots \\[1mm]
\, \epsilon_0^{n-1} \, & \cdots & \epsilon_{n-1}^{n-1} & \, g_{n-1}^{} \,
\end{pmatrix}
\begin{pmatrix}
\omega_0^{} \\[1mm]
\vdots \\[1mm]
\, \omega_{n-1}^{} \, \\[1mm]
\, \omega_{n}^{} \,
\end{pmatrix}
~=~
\begin{pmatrix}
0 \\[1mm]
h_1^{} \\[1mm]
\vdots \\[1mm]
%\, 0 \,
\, h_n^{} \, 
\end{pmatrix}
\label{eq:generalConstraints}
\end{align}

\noindent
If $g_{j}^{} = (\alpha\, \epsilon_n^{})^j$ for some free variable $\alpha$ and $h_{j}^{}=0$, it follows from \cref{eq:limitSolution} that the parameter set below has asymptotic behaviours satisfying the above constraints.
\begin{align}
\omega_{i (\neq n)}^{} ~=~
\frac{1}{f_i^{} - \alpha f_n^{}} \prod_{j \neq \{ i, n \}}^{} \frac{1}{f_i^{} - f_j^{} } 
\hspace{12mm}
\omega_n^{} ~=~
%(-1)^{n}\, 
\prod_{j \neq n}^{} \, \frac{1}{\alpha f_n^{} - f_j^{}} 
\label{eq:alternative}
\end{align}

\noindent
The limit of coinciding pairs of interpolation nodes could also be considered in order to find parameters suitable for the iteration scheme with first derivatives. However, this paper does not seek to comprehensively explore such cases. It is simply emphasised that other parameter sets can be applied that have similar convergence properties.

%not as convenient to update

\vfill\pagebreak
%\vspace{5mm}\phantom{.}

\begin{table}[!htbp]
	\centering
	%\tiny
	%\scriptsize
	%\footnotesize
	\small
	%\normalsize
	\begin{tabular}{llp{0mm}lp{0mm}l} %p{18mm}
		\toprule
		\\[-2mm]
		\multicolumn{2}{l}{ Method: }
		&& \multicolumn{3}{l}{ $\displaystyle x_{n+1}^{} ~=\,
			\left.
			\left( \, \sum_{i=0}^{n} \, \frac{\omega_i^{} \, x_i^{}}{f_i^{}} \, \right) \! \middle/ \! \left( \, \sum_{i=0}^{n} \, \frac{\omega_i^{}}{f_i^{}} \, \right) \right.$ }
		\\[7mm] %\midrule %\midrule
		%\\[-2mm]
		\multicolumn{2}{l}{ Weights: }
		&& $\displaystyle \omega_i^{} ~=~ \prod_{j\neq i}^{} \, \frac{1}{x_i^{} - x_j^{}}$
		&& $\displaystyle \omega_i^{} ~=~ \prod_{j\neq i}^{} \, \frac{1}{f_i^{} - f_j^{}}$
		\\[7mm]
		\multicolumn{2}{l}{ Interpolant:}
		&& 
		$f[x]$: (1,n-1) rational 
		%\begin{tabular}{@{}ll}
		%	$f[x]$: (1,n-1) rational \\
		%	$x[f]$: (n,n) rational \\
		%\end{tabular}
		&& 
		$f[x]$: (n+1,n-1) rational
		%\begin{tabular}{@{}ll}
		%	$f[x]$: (n+1,n-1) rational \\
		%	$x[f]$: (n) polynomial \\
		%\end{tabular}
		\\[1mm]
		&
		&& $x[f]$: (n,n) rational
		&& $x[f]$: (n) polynomial 
		\\[2mm] \midrule \midrule
		n+1 &  
		&& Leading error factor $(E_{1,n+1}^{})$ && Leading error factor $(E_{1,n+1}^{})$ \\ \midrule
		\\[-2mm]
		2
		& %$\displaystyle \left( f_\ast^{(2)} \right)$
		&& $\displaystyle \frac{f_\ast^{(2)}}{2 \, f_\ast^{(1)}} $
		&& $\displaystyle \frac{f_\ast^{(2)}}{2 \, f_\ast^{(1)}}  $
		\\[7mm]
		3
		& %$\displaystyle \left( f_\ast^{(2,2)}, f_\ast^{(1,3)} \right)$
		&& $\displaystyle \frac{3\, f_\ast^{(2,2)} - 2\, f_\ast^{(1,3)}}{12 \, f_\ast^{(1,1)}} $
		&& $\displaystyle \frac{6\, f_\ast^{(2,2)} - 2\, f_\ast^{(1,3)}}{12 \, f_\ast^{(1,1)}} $
		\\[7mm]
		4
		& %$\displaystyle \left( f_\ast^{(2,2,2)}, f_\ast^{(1,2,3)}, f_\ast^{(1,1,4)} \right)$
		&& $\displaystyle \frac{3\, f_\ast^{(2,2,2)} - 4\, f_\ast^{(1,2,3)}  + f_\ast^{(1,1,4)}}{24 \, f_\ast^{(1,1,1)}} $
		&& $\displaystyle \frac{15\, f_\ast^{(2,2,2)} - 10\, f_\ast^{(1,2,3)}  + f_\ast^{(1,1,4)}}{24 \, f_\ast^{(1,1,1)}} $
		\\[5mm] \bottomrule
	\end{tabular}
	\vspace{3mm} %\small 
	\caption{Root search schemes based on derivative-free rational interpolation.}
		%(error notation defined in eq~\ref{eq:error} and \ref{eq:notationfi}).}
	\label{tbl:barycentric0}
\end{table}

\vfill

\begin{table}[!htbp]
	\centering
	%\tiny
	%\scriptsize
	%\footnotesize
	\small
	%\normalsize
	\begin{tabular}{p{3mm}p{20mm}p{20mm}p{20mm}p{20mm}p{20mm}}
		\toprule
		& Picard & Secant & & & Newton \\[0mm] %\midrule
		$i$ &  & $n=1$ & $n=2$ & $n=3$ &  \\ \midrule
$ 0 $ & 	$ 2.26 $ & 	$ 2.26 $ & 	$ 2.26 $ & 	$ 2.26 $ & 	$ 2.26 $ \\
$ 1 $ & 	$ 1.73 $ & 	$ 1.73 $ & 	$ 1.73 $ & 	$ 1.73 $ & 	$ 1.24 $ \\
$ 2 $ & 	$ 1.90\!\times\! 10^{-1} $ & 	$ 6.19\!\times\! 10^{-1} $ & 	$ 6.19\!\times\! 10^{-1} $ & 	$ 6.19\!\times\! 10^{-1} $ & 	$ 1.39 $ \\
$ 3 $ & 	$ 1.14\!\times\! 10^{-1} $ & 	$ 8.35\!\times\! 10^{-1} $ & 	$ 3.47\!\times\! 10^{-1} $ & 	$ 3.47\!\times\! 10^{-1} $ & 	$ 4.94\!\times\! 10^{-2} $ \\
$ 4 $ & 	$ 8.15\!\times\! 10^{-2} $ & 	$ 1.01\!\times\! 10^{-1} $ & 	$ 6.61\!\times\! 10^{-2} $ & 	$ 1.77\!\times\! 10^{-2} $ & 	$ 5.68\!\times\! 10^{-4} $ \\
$ 5 $ & 	$ 5.24\!\times\! 10^{-2} $ & 	$ 1.23\!\times\! 10^{-2} $ & 	$ 1.73\!\times\! 10^{-3} $ & 	$ 2.00\!\times\! 10^{-4} $ & 	$ 7.12\!\times\! 10^{-8} $ \\
$ 6 $ & 	$ 3.63\!\times\! 10^{-2} $ & 	$ 2.91\!\times\! 10^{-4} $ & 	$ 4.27\!\times\! 10^{-6} $ & 	$ 1.78\!\times\! 10^{-8} $ & 	$ 1.12\!\times\! 10^{-15} $ \\
$ 7 $ & 	$ 2.40\!\times\! 10^{-2} $ & 	$ 7.94\!\times\! 10^{-7} $ & 	$ 5.60\!\times\! 10^{-11} $ & 	$ 4.40\!\times\! 10^{-16} $ & 	$ 2.76\!\times\! 10^{-31} $ \\
$ 8 $ & 	$ 1.63\!\times\! 10^{-2} $ & 	$ 5.09\!\times\! 10^{-11} $ & 	$ 4.80\!\times\! 10^{-20} $ & 	$ 6.06\!\times\! 10^{-31} $ & 	$ 1.68\!\times\! 10^{-62} $ \\
$ 9 $ & 	$ 1.09\!\times\! 10^{-2} $ & 	$ 8.93\!\times\! 10^{-18} $ & 	$ 1.33\!\times\! 10^{-36} $ & 	$ 2.08\!\times\! 10^{-59} $ & 	$ 6.25\!\times\! 10^{-125} $ \\
%$ 10 $ & 	$ 7.39\!\times\! 10^{-3} $ & 	$ 1.00\!\times\! 10^{-28} $ & 	$ 4.14\!\times\! 10^{-67} $ & 	$ 2.17\!\times\! 10^{-114} $ & 	$ 8.63\!\times\! 10^{-250} $ \\
		\bottomrule
	\end{tabular}
	\vspace{2mm}
	\caption{Error magnitudes for the root of $f[x] = \cos x - x$ ($x_0^{} = 3$), where
		Picard's method corresponds to fixed-point iteration: $x_{i+1}^{} = \cos x_i^{}$. %\newline
		Newton's method is applied in the right column using $f_i^{\prime}$ information,
		and the central sequences are generated by the method and weights defined in the left column of \cref{tbl:barycentric0} using only the latest $n+1$ points.
	}
	\label{tbl:example0}
\end{table}

\vfill\pagebreak
%\vspace{5mm}\phantom{.}

\begin{table}[!htbp]
	\centering
	%\tiny
	%\scriptsize
	%\footnotesize
	\small
	%\normalsize
	\begin{tabular}{llp{0mm}lp{0mm}l} %p{18mm}
		\toprule
		\\[-2mm]
		\multicolumn{2}{l}{ Method: }
		&& \multicolumn{3}{l}{ $\displaystyle x_{n+1}^{} ~=\,
			\left.
			\left( \, \sum_{i=0}^{n} \, \frac{ \lambda_i^{} \, ( x_i^{}  - f_i^{} / f_i^\prime ) - \gamma_i^{} \, f_i^{} \, x_i^{} }{ f_i^2 } \, \right) \! \middle/ \! \left( \, \sum_{i=0}^{n} \, \frac{ \lambda_i^{}  - \gamma_i^{} \, f_i^{} }{f_i^2} \, \right) \right.$ }
		\\[7mm] %\midrule %\midrule
		%\\[-2mm]
		\multicolumn{2}{l}{ Weights: }
		&& $\displaystyle \lambda_i^{} ~=~ f_i^{\prime} ~ \prod_{j\neq i}^{} \, \frac{1}{(x_i^{} - x_j^{})^2}$
		&& $\displaystyle \lambda_i^{} ~=~ \prod_{j\neq i}^{} \, \frac{1}{(f_i^{} - f_j^{})^2}$
		\\[7mm]
		\multicolumn{2}{l}{  }
		&& $\displaystyle \gamma_i^{} ~=~ - \frac{2\, \lambda_i^{}}{f_i^{\prime}} \, \sum_{j\neq i}^{} \, \frac{1}{x_i^{} - x_j^{}}$
		&& $\displaystyle \gamma_i^{} ~=~ - 2\, \lambda_i^{}\, \sum_{j\neq i}^{} \, \frac{1}{f_i^{} - f_j^{}}$
		\\[7mm]
		\multicolumn{2}{l}{ Interpolant:}
		&& $f[x]$: (1,2n) rational 
		&& $f[x]$: (2n+2,2n) rational
		\\[1mm]
		&
		&& $x[f]$: (2n+1,2n+1) rational
		&& $x[f]$: (2n+1) polynomial 
		\\[2mm] \midrule \midrule
		n+1 &  
		&& Leading error factor $(E_{2,n+1}^{})$ && Leading error factor $(E_{2,n+1}^{})$ \\ \midrule
		\\[-2mm]
		1
		& %$\displaystyle \left( f_\ast^{(2)} \right)$
		&& $\displaystyle \frac{f_\ast^{(2)}}{2 \, f_\ast^{(1)}} $
		&& $\displaystyle \frac{f_\ast^{(2)}}{2 \, f_\ast^{(1)}}  $
		\\[7mm]
		2
		& %$\displaystyle \left( f_\ast^{(2,2,2)}, f_\ast^{(1,2,3)}, f_\ast^{(1,1,4)} \right)$
		&& $\displaystyle \frac{3\, f_\ast^{(2,2,2)} - 4\, f_\ast^{(1,2,3)}  + f_\ast^{(1,1,4)}}{24 \, f_\ast^{(1,1,1)}} $
		&& $\displaystyle \frac{15\, f_\ast^{(2,2,2)} - 10\, f_\ast^{(1,2,3)}  + f_\ast^{(1,1,4)}}{24 \, f_\ast^{(1,1,1)}} $
		\\[5mm] \bottomrule
	\end{tabular}
	\vspace{3mm} %\small 
	\caption{Root search schemes based on rational interpolation with first~derivatives.}
	%(error notation defined in eq~\ref{eq:error} and \ref{eq:notationfi}).}
	\label{tbl:barycentric1}
\end{table}

\vfill

\begin{table}[!htbp]
	\centering
	%\tiny
	%\scriptsize
	%\footnotesize
	\small
	%\normalsize
	\begin{tabular}{p{3mm}p{20mm}p{20mm}p{20mm}p{20mm}p{20mm}}
		\toprule
		& Newton & & & & Halley \\[0mm] %\midrule
		$i$ & $n=0$ & $n=1$ & $n=2$ & $n=3$ &  \\ \midrule
$ 0 $ & 	$ 2.26 $ & 	$ 2.26 $ & 	$ 2.26 $ & 	$ 2.26 $ & 	$ 2.26 $ \\
$ 1 $ & 	$ 1.24 $ & 	$ 1.24 $ & 	$ 1.24 $ & 	$ 1.24 $ & 	$ 8.72\!\times\! 10^{-1} $ \\
$ 2 $ & 	$ 1.39 $ & 	$ 1.18\!\times\! 10^{-1} $ & 	$ 1.18\!\times\! 10^{-1} $ & 	$ 1.18\!\times\! 10^{-1} $ & 	$ 5.27\!\times\! 10^{-2} $ \\
$ 3 $ & 	$ 4.94\!\times\! 10^{-2} $ & 	$ 6.85\!\times\! 10^{-4} $ & 	$ 2.44\!\times\! 10^{-5} $ & 	$ 2.44\!\times\! 10^{-5} $ & 	$ 1.65\!\times\! 10^{-5} $ \\
$ 4 $ & 	$ 5.68\!\times\! 10^{-4} $ & 	$ 1.35\!\times\! 10^{-10} $ & 	$ 9.33\!\times\! 10^{-15} $ & 	$ 4.76\!\times\! 10^{-15} $ & 	$ 5.19\!\times\! 10^{-16} $ \\
$ 5 $ & 	$ 7.12\!\times\! 10^{-8} $ & 	$ 1.88\!\times\! 10^{-28} $ & 	$ 2.87\!\times\! 10^{-43} $ & 	$ 6.73\!\times\! 10^{-44} $ & 	$ 1.62\!\times\! 10^{-47} $ \\
$ 6 $ & 	$ 1.12\!\times\! 10^{-15} $ & 	$ 1.41\!\times\! 10^{-77} $ & 	$ 1.56\!\times\! 10^{-126} $ & 	$ 7.76\!\times\! 10^{-131} $ & 	$ 4.93\!\times\! 10^{-142} $ \\
%$ 7 $ & 	$ 2.76\!\times\! 10^{-31} $ & 	$ 1.55\!\times\! 10^{-211} $ & 	$ 6.41\!\times\! 10^{-369} $ & 	$ 3.67\!\times\! 10^{-389} $ & 	$ 1.39\!\times\! 10^{-425} $ \\
		\bottomrule
	\end{tabular}
	\vspace{2mm}
	\caption{Error magnitudes for the root of $f[x] = \cos x - x$ ($x_0^{} = 3$). Halley's method is applied in the right column using $f_i^{\prime\prime}$ information, and the other sequences are generated by the method and weights defined in the left column of \cref{tbl:barycentric1} using only the latest $n+1$ points.
	}
\label{tbl:example1}
\end{table}

\vfill\pagebreak
\section{Approximate roots of univariate rational interpolants}

In the previous section, the exact roots of interpolating functions were selected for subsequent iteration. However, such functions fundamentally are approximations. An alternative approach is to apply local iteration schemes about a specific point of the interpolant. The local scheme then approximates the interpolant root, but if the associated errors are of similar order to the interpolation errors, the iteration should be similarly effective.

\Cref{tbl:interp} summarises the choices for derivative-free interpolation made in this section. The derivatives are deduced by performing a series expansion about the respective points. The leading error terms for Newton iteration ($x_{n+1}^{} = x_{n}^{} - f_n^{}/f_n^\prime$) with interpolant derivatives are then also presented for arbitrary~$\omega_i^{}$. 
%
%\vspace{1mm}
\begin{table}[!htbp]
	\centering
	%\tiny
	%\scriptsize
	%\footnotesize
	\small
	%\normalsize
	\begin{tabular}{l@{\hspace{4mm}}l} %p{18mm}
		\\[0mm]
		\toprule
		Direct function interpolation & Inverse function interpolation
		\\ \midrule
		\\[-2mm]
		$
		\begin{array}{@{}l@{\hspace{2mm}}l@{\hspace{1mm}}l@{}}
		\displaystyle
		f[x] &\displaystyle \approx& \, \displaystyle
		\left.
		\left( \, \sum_{i=0}^{n} \, \frac{\omega_i^{} \, f_i^{}}{x-x_i^{}} \, \right) \! \middle/ \! \left( \, \sum_{i=0}^{n} \, \frac{\omega_i^{}}{x-x_i^{}} \, \right) \right.
		\\[7mm]
		\displaystyle
		f_i^\prime &\displaystyle \approx& \displaystyle
		~- \frac{1}{\omega_i^{}} \, \sum_{j \neq i}^{} \, \omega_j^{} \, \frac{ f_i^{}-f_j^{} }{x_i^{} - x_j^{} } 
		\end{array}$
		&
		$
		\begin{array}{@{}l@{\hspace{2mm}}l@{\hspace{1mm}}l@{}}
		\displaystyle
		x[f] &\displaystyle \approx& \, \displaystyle
		\left.
		\left( \, \sum_{i=0}^{n} \, \frac{\omega_i^{} \, x_i^{}}{f-f_i^{}} \, \right) \! \middle/ \! \left( \, \sum_{i=0}^{n} \, \frac{\omega_i^{}}{f-f_i^{}} \, \right) \right.
		\\[7mm]
		\displaystyle
		\frac{1}{ f_i^\prime }  &\displaystyle \approx& \displaystyle
		~- \frac{1}{\omega_i^{}} \, \sum_{j \neq i}^{} \, \omega_j^{} \, \frac{ x_i^{}-x_j^{} }{f_i^{} - f_j^{} }
		\end{array}$
		\\[11mm] \midrule \midrule
		\\[-2mm]
		$\displaystyle
		\epsilon_{n+1}^{} ~\sim~
		\epsilon_n^{} 
		\left.
		\Bigg( \, \sum_{i=0}^{n} \, \omega_i^{}  \, \Bigg) \! \middle/ \! \Bigg( \, \sum_{j \neq n}^{} \, \omega_j^{}  \, \Bigg) \right.$
		%\, +\, \cdots \text{higher-order terms}$
		&
		$\displaystyle
		\epsilon_{n+1}^{} ~\sim~
		\Bigg( \, \sum_{i=0}^{n} \, \omega_i^{} \, \Bigg) \, \frac{\epsilon_n^{} }{\omega_n^{}}$
		%\, +\, \cdots \text{higher-order terms}$
		\\[5mm] \bottomrule
	\end{tabular}
	\vspace{3mm} %\small 
	\caption{Interpolants with corresponding first derivatives, and leading error terms for Newton iteration.}
	\label{tbl:interp} %about $x_n^{}
\end{table}

%\vspace{-3mm}
\noindent
As in the previous section, the leading error is generally $\mathcal{O}[\epsilon]$, but that term vanishes when $\sum_{i}^{} \, \omega_i^{} = 0$. In such cases, the interpolant derivative becomes a weighted average of finite differences. 
If then successively requiring the leading error to vanish up to $\mathcal{O}[\epsilon_n^{2}]$ terms, the following constraints also apply:
\begin{align}
\sum_{j\neq n}^{} \, \omega_j^{} \, \epsilon_j^{k} ~=~ \mathcal{O}[\epsilon_n^{}]
%0
\hspace{15mm} 
%0 \leq m \leq n-1
1 \leqslant k\leqslant n-1
\label{eq:constraintsA}
\end{align}

\noindent
The constraints take the same form as \cref{eq:generalConstraints}, and so the weights identified in the previous section can again be used. Similarly, it follows that when these constraints are satisfied, the leading error term takes the form of \cref{eq:error} for which the asymptotic order of convergence was previously discussed.

%asymptotic order of convergence is as stated in table~\ref{tbl:rootConvergence}.

\Cref{tbl:Xapprox0,tbl:Fapprox0} present iteration schemes based on the above interpolating functions. The leading error factors for the schemes in \cref{tbl:Xapprox0} are common to those when taking the exact interpolant root of $x[f]$ (\cref{tbl:barycentric0}). However, the sub-dominant error terms differ, so the iteration sequences are generally distinct. 
The scheme in the left column of \cref{tbl:Fapprox0} is also noted to be equivalent to the method described by Sidi~\cite{sidi2008generalization}, though with the polynomial interpolant derivative being constructed here without recursion. This method should be generally favoured amongst the derivative-free root search methods detailed in this paper, as its leading error factor contains only the highest degree term common to all methods, and so the error tends to be suppressed.

%recursively calculating divided differences

\vfill\pagebreak
%\vspace{0mm}\phantom{.}

\begin{table}[!htbp]
	\centering
	%\tiny
	%\scriptsize
	%\footnotesize
	\small
	%\normalsize
	\begin{tabular}{llp{0mm}lp{0mm}l} %p{18mm}
		\toprule
		\\[-2mm]
		\multicolumn{2}{l}{ Method: }
		&& \multicolumn{3}{l}{ $\displaystyle x_{n+1}^{} ~=~ 
			%x_{n}^{} - \alpha_n^{} \, f_{n}^{} $ \hspace{2mm} : \hspace{1mm}
			x_{n}^{} - \frac{ f_{n}^{} }{ f_n^\prime } $ \hspace{2mm} : \hspace{1mm}
			$\displaystyle 
			%\alpha_n^{} ~=~ 
			\frac{ 1 }{ f_n^\prime } ~\approx\,
			\left.
			\Bigg( \, \sum_{k\neq n}^{} \, \omega_k^{}  \, \frac{ x_n^{} - x_k^{} }{ f_n^{} - f_k^{} } \, \Bigg) \! \middle/ \! \Bigg( \, \sum_{k\neq n}^{} \, \omega_k^{} \, \Bigg) \right.$
		}
		\\[7mm] %\midrule %\midrule
		%\\[-2mm]
		\multicolumn{2}{l}{ Weights: }
		&& %$\displaystyle \omega_j^{} ~=~ \frac{1}{f_j^{}}  \sum_{k\neq \{j,n\} }^{}  \frac{1}{x_j^{} - x_k^{}}$
		$\displaystyle \omega_i^{} ~=~ \prod_{j\neq i}^{} \, \frac{1}{x_i^{} - x_j^{}}$
		&& $\displaystyle \omega_i^{} ~=~ \prod_{j\neq i}^{} \, \frac{1}{f_i^{} - f_j^{}}$
		\\[7mm]
		\multicolumn{2}{l}{ Interpolant:}
		&& 
		$x[f]$: (n,n) rational 
		%\begin{tabular}{@{}ll}
		%	$f[x]$: (1,n-1) rational \\
		%	$x[f]$: (n,n) rational \\
		%\end{tabular}
		&& 
		$x[f]$: (n) polynomial
		%\begin{tabular}{@{}ll}
		%	$f[x]$: (n+1,n-1) rational \\
		%	$x[f]$: (n) polynomial \\
		%\end{tabular}
		%%%%%%%%%%%
		%\\[1mm]
		%&
		%&& $x[f]$: (n,n) rational
		%&& $x[f]$: (n) polynomial 
		\\[2mm] \midrule \midrule
		n+1 &  
		&& Leading error factor $(E_{1,n+1}^{})$ && Leading error factor $(E_{1,n+1}^{})$ \\ \midrule
		\\[-2mm]
		2
		& %$\displaystyle \left( f_\ast^{(2)} \right)$
		&& $\displaystyle \frac{f_\ast^{(2)}}{2 \, f_\ast^{(1)}} $
		&& $\displaystyle \frac{f_\ast^{(2)}}{2 \, f_\ast^{(1)}}  $
		\\[7mm]
		3
		& %$\displaystyle \left( f_\ast^{(2,2)}, f_\ast^{(1,3)} \right)$
		&& $\displaystyle \frac{3\, f_\ast^{(2,2)} - 2\, f_\ast^{(1,3)}}{12 \, f_\ast^{(1,1)}} $
		&& $\displaystyle \frac{6\, f_\ast^{(2,2)} - 2\, f_\ast^{(1,3)}}{12 \, f_\ast^{(1,1)}} $
		\\[7mm]
		4
		& %$\displaystyle \left( f_\ast^{(2,2,2)}, f_\ast^{(1,2,3)}, f_\ast^{(1,1,4)} \right)$
		&& $\displaystyle \frac{3\, f_\ast^{(2,2,2)} - 4\, f_\ast^{(1,2,3)}  + f_\ast^{(1,1,4)}}{24 \, f_\ast^{(1,1,1)}} $
		&& $\displaystyle \frac{15\, f_\ast^{(2,2,2)} - 10\, f_\ast^{(1,2,3)}  + f_\ast^{(1,1,4)}}{24 \, f_\ast^{(1,1,1)}} $
		\\[5mm] \bottomrule
	\end{tabular}
	\vspace{3mm} %\small 
	\caption{Root search schemes based on Newton's method for $x[f]$ interpolant.}
	%(error notation defined in eq~\ref{eq:error} and \ref{eq:notationfi}).}
	%based on derivative-free interpolation of the inverse function $x[f]$.}
	\label{tbl:Xapprox0}
\end{table}

\vfill
%\vspace{-2mm}
%\pagebreak
\begin{table}[!htbp]
	\centering
	%\tiny
	%\scriptsize
	%\footnotesize
	\small
	%\normalsize
	%\vspace{-10mm}
	\begin{tabular}{llp{0mm}lp{0mm}l} %p{18mm}
		\toprule
		\\[-2mm]
		\multicolumn{2}{l}{ Method: }
		&& \multicolumn{3}{l}{ $\displaystyle x_{n+1}^{} ~=~ 
			%x_{n}^{} - \alpha_n^{} \, f_{n}^{} $ \hspace{2mm} : \hspace{1mm}
			x_{n}^{} - \frac{ f_{n}^{} }{ f_n^\prime } $ \hspace{2mm} : \hspace{1mm}
			$\displaystyle 
			%\alpha_n^{-1} ~=~ 
			f_n^\prime ~\approx\,
			\left.
			\Bigg( \, \sum_{k\neq n}^{} \, \omega_k^{}  \, \frac{ f_n^{} - f_k^{} }{ x_n^{} - x_k^{} } \, \Bigg) \! \middle/ \! \Bigg( \, \sum_{k\neq n}^{} \, \omega_k^{} \, \Bigg) \right.$
		}
		\\[7mm] %\midrule %\midrule
		%\\[-2mm]
		\multicolumn{2}{l}{ Weights: }
		&& %$\displaystyle \omega_j^{} ~=~ \frac{1}{f_j^{}}  \sum_{k\neq \{j,n\} }^{}  \frac{1}{x_j^{} - x_k^{}}$
		$\displaystyle \omega_i^{} ~=~ \prod_{j\neq i}^{} \, \frac{1}{x_i^{} - x_j^{}}$
		&& $\displaystyle \omega_i^{} ~=~ \prod_{j\neq i}^{} \, \frac{1}{f_i^{} - f_j^{}}$
		\\[7mm]
		\multicolumn{2}{l}{ Interpolant:}
		&& 
		$f[x]$: (n) polynomial
		%\begin{tabular}{@{}ll}
		%	$f[x]$: (1,n-1) rational \\
		%	$x[f]$: (n,n) rational \\
		%\end{tabular}
		&& 
		$f[x]$: (n,n) rational 
		%\begin{tabular}{@{}ll}
		%	$f[x]$: (n+1,n-1) rational \\
		%	$x[f]$: (n) polynomial \\
		%\end{tabular}
		%%%%%%%%%%%
		%\\[1mm]
		%&
		%&& $x[f]$: (n,n) rational
		%&& $x[f]$: (n) polynomial 
		\\[2mm] \midrule \midrule
		n+1 &  
		&& Leading error factor $(E_{1,n+1}^{})$ && Leading error factor $(E_{1,n+1}^{})$ \\ \midrule
		\\[-2mm]
		2
		& %$\displaystyle \left( f_\ast^{(2)} \right)$
		&& $\displaystyle \frac{f_\ast^{(2)}}{2 \, f_\ast^{(1)}} $
		&& $\displaystyle \frac{f_\ast^{(2)}}{2 \, f_\ast^{(1)}}  $
		\\[7mm]
		3
		& %$\displaystyle \left( f_\ast^{(2,2)}, f_\ast^{(1,3)} \right)$
		&& $\displaystyle \frac{-\, f_\ast^{(3)}}{6 \, f_\ast^{(1)}} $
		&& $\displaystyle \frac{3\, f_\ast^{(2,2)} - 2\, f_\ast^{(1,3)}}{12 \, f_\ast^{(1,1)}} $
		\\[7mm]
		4
		& %$\displaystyle \left( f_\ast^{(2,2,2)}, f_\ast^{(1,2,3)}, f_\ast^{(1,1,4)} \right)$
		&& $\displaystyle \frac{f_\ast^{(4)}}{24 \, f_\ast^{(1)}} $
		&& $\displaystyle \frac{6\, f_\ast^{(2,2,2)} - 6\, f_\ast^{(1,2,3)}  + f_\ast^{(1,1,4)}}{24 \, f_\ast^{(1,1,1)}} $
		\\[5mm] \bottomrule
	\end{tabular}
	\vspace{3mm} %\small 
	\caption{Root search schemes based on Newton's method for $f[x]$ interpolant.}
		%(error notation defined in eq~\ref{eq:error} and \ref{eq:notationfi}).}
		%based on derivative-free interpolation of the direct function $f[x]$.}
	\label{tbl:Fapprox0}
\end{table}

\vfill\pagebreak

%\vfill
%\vspace{5mm}\phantom{.}
\begin{table}[!htbp]
	\centering
	%\tiny
	%\scriptsize
	%\footnotesize
	\small
	%\normalsize
	\begin{tabular}{@{}llp{0mm}lp{1mm}l@{}} %p{18mm}
		\toprule
		\\[-2mm]
		\multicolumn{2}{l}{ Method: }
		&& \multicolumn{3}{l}{ $\displaystyle x_{n+1}^{} ~=~ 
			%x_{n}^{} - \alpha_n^{} \, \phi_{n}^{\prime} $ %\hspace{2mm} : \hspace{1mm}
			x_n^{}
			- \left( \frac
			{ f_n^{\prime \,2} + (\frac12 - \beta ) \, f_n^{} \, f_n^{\prime\prime} }
			{ f_n^{\prime \,2} - \beta \, f_n^{} \, f_n^{\prime\prime}  } \right) 
			\frac{ f_n^{} }{ f_n^{\prime} }
			%\hspace{12mm} %: \hspace{1mm} 
			\hfill
			\begin{array}{@{}l@{}}
			\beta ~\text{free parameter} \\[0.5mm]
			\beta = 1 ~\text{recommended}
			\end{array}
			$
		}
		\\[8mm] %\midrule %\midrule
		&&& \multicolumn{3}{l}{ %\hspace{1mm}
			$\displaystyle 
			f_i^{\prime \prime}
			~\approx~
			\frac{2  f_i^{\prime \, 3}}{\lambda_i^{}}
			\left( \frac{\gamma_i^{} }{ f_i^{\prime} }
			+\sum_{k \neq i}^{} \frac
			{ \lambda_k^{} \, (x_i^{} - x_k^{}) + ( \gamma_k^{} \, (x_i^{} - x_k^{}) - \lambda_k^{} / f_k^{\prime}) \, (f_i^{} - f_k^{}) }
			{ (f_i^{} - f_k^{})^2 }
			\right)
			$ }
		\\[8mm]
		%\\[-2mm]
		\multicolumn{2}{l}{ Weights: }
		&& $\displaystyle \lambda_i^{} ~=~ \prod_{j\neq i}^{} \, \frac{1}{(f_i^{} - f_j^{})^2}$
		&& %\hspace{1mm}
		$\displaystyle \gamma_i^{} ~=~ -\, 2\, \lambda_i^{} \, \sum_{j\neq i}^{} \, \frac{1}{f_i^{} - f_j^{}}$
		\\[5mm] \midrule \midrule
		\\[-3mm]
		\multicolumn{2}{l}{ Error: }
		&&
		\multicolumn{2}{l}{% Leading errors: \hspace{4mm} 
			$\displaystyle \epsilon_{1}^{} ~\sim~ 
			\frac{f_\ast^{(2)}}{2 \, f_\ast^{(1)}} \, \epsilon_{0}^{2} $
			%\hspace{8mm} $\displaystyle \epsilon_{3}^{} ~\sim~ \frac{f_\ast^{(4)}}{24 \, f_\ast^{(1)}}\, \epsilon_0^{2}\, \epsilon_1^{2} $ 
			%\hfill ($\beta$-independent)
			\hspace{5mm}\phantom{.}
		}
		&
		\multicolumn{1}{l}{
			$\displaystyle \epsilon_{2}^{} ~\sim~ \frac{15\, f_\ast^{(2,2,2)} - 10\, f_\ast^{(1,2,3)}+ f_\ast^{(1,1,4)}}{24 \, f_\ast^{(1,1,1)}}\, \epsilon_0^{2}\, \epsilon_1^{2} $ 
			%\hfill ($\beta$-independent)
		}
		\\[4mm] \bottomrule
	\end{tabular}
	\vspace{3mm} %\small %$f_\ast^{\,\cdots}$
	\caption{Root search scheme based on Chebyshev-Halley methods for $x[f]$ interpolant.}
		%(error notation defined in eq~\ref{eq:notationfi}).}
	\label{tbl:Xapprox1}
\end{table}

%\vfill

\begin{table}[!htbp]
	\centering
	%\tiny
	%\scriptsize
	%\footnotesize
	\small
	%\normalsize
	\begin{tabular}{@{}llp{0mm}lp{1mm}l@{}} %p{18mm}
		\toprule
		\\[-2mm]
		\multicolumn{2}{l}{ Method: }
		&& \multicolumn{3}{l}{ $\displaystyle x_{n+1}^{} ~=~ 
			%x_{n}^{} - \alpha_n^{} \, \phi_{n}^{\prime} $ %\hspace{2mm} : \hspace{1mm}
			x_n^{}
			- \left( \frac
			{ f_n^{\prime \,2} + (\frac12 - \beta ) \, f_n^{} \, f_n^{\prime\prime} }
			{ f_n^{\prime \,2} - \beta \, f_n^{} \, f_n^{\prime\prime}  } \right) 
			\frac{ f_n^{} }{ f_n^{\prime} }
			%\hspace{12mm} %: \hspace{1mm} 
			\hfill
			\begin{array}{@{}l@{}}
			\beta ~\text{free parameter} \\[0.5mm]
			\beta = 1 ~\text{recommended}
			\end{array}
			$
		}
		\\[8mm] %\midrule %\midrule
		&&& \multicolumn{3}{l}{ %\hspace{1mm}
			$\displaystyle 
			f_i^{\prime \prime}
			~\approx~
			- \frac{2}{\lambda_i^{}}
			\left( \gamma_i^{} \,f_i^{\prime} 
			+\sum_{k \neq i}^{} \frac
			{ \lambda_k^{} \, (f_i^{} - f_k^{}) + ( \gamma_k^{} \, (f_i^{} - f_k^{}) - \lambda_k^{}\, f_k^{\prime}) \, (x_i^{} - x_k^{}) }
			{ (x_i^{} - x_k^{})^2 }
			\right)
			$ }
		\\[8mm]
		%\\[-2mm]
		\multicolumn{2}{l}{ Weights: }
		&& $\displaystyle \lambda_i^{} ~=~ \prod_{j\neq i}^{} \, \frac{1}{(x_i^{} - x_j^{})^2}$
		&& %\hspace{1mm}
		$\displaystyle \gamma_i^{} ~=~ - 2\, \lambda_i^{} \, \sum_{j\neq i}^{} \, \frac{1}{x_i^{} - x_j^{}}$
		\\[5mm] \midrule \midrule
		\\[-3mm]
		\multicolumn{2}{l}{ Error: }
		&&
		\multicolumn{2}{l}{%Leading errors: \hspace{4mm} 
			$\displaystyle \epsilon_{1}^{} ~\sim~ 
			\frac{f_\ast^{(2)}}{2 \, f_\ast^{(1)}} \, \epsilon_{0}^{2} $
			%\hspace{8mm} $\displaystyle \epsilon_{3}^{} ~\sim~ \frac{f_\ast^{(4)}}{24 \, f_\ast^{(1)}}\, \epsilon_0^{2}\, \epsilon_1^{2} $ 
			%\hfill ($\beta$-independent)
			\hspace{5mm}\phantom{.}
		}
		&
		\multicolumn{1}{l}{
			$\displaystyle \epsilon_{2}^{} ~\sim~ \frac{f_\ast^{(4)}}{24 \, f_\ast^{(1)}}\, \epsilon_0^{2}\, \epsilon_1^{2} $ 
			\hfill ($\beta$-independent) \hspace{5mm}
		}
		\\[4mm] \bottomrule
	\end{tabular}
	\vspace{3mm} %\small 
	\caption{Root search scheme based on Chebyshev-Halley methods for $f[x]$ interpolant.}
		%(error notation defined in eq~\ref{eq:notationfi}).}
	\label{tbl:Fapprox1}
\end{table}

\Cref{tbl:Xapprox1,tbl:Fapprox1} present iteration schemes based on interpolation with first derivatives, using \cref{eq:xfHermite} for $x[f]$ and its counterpart form for $f[x]$. The leading error factors in \cref{tbl:Xapprox1} are common to those when taking the exact interpolant root of $x[f]$ (\cref{tbl:barycentric0}). However, the iteration sequences are generally distinct. 
In order to deduce the second derivative $f_i^{\prime\prime}$ from the $x[f]$ interpolant, $x^{\prime\prime}[f_i^{}]$ can be first determined by a series expansion about $f_i^{}$ and then the following relation used:
\begin{align}
\frac{\partial^2 x}{ \partial f^2}
~=~
\frac{\partial}{\partial f} \left( \frac{\partial x}{\partial f} \right)
\,=~
\frac{\partial x}{\partial f} \frac{\partial }{\partial x} \left( \frac{1}{f^\prime} \right)
\,=\,
- \frac{f^{\prime\prime}}{f^{\prime \, 3}}
\end{align}

\noindent
The Chebyshev-Halley family~\cite{gutierrez1997family} of root search methods is applied to the interpolant, as these are known to be cubically convergent when the true second derivative is used. For the iteration schemes based on exact interpolant roots, the orders of convergence were identified to be super-quadratic (for $n >0$) but sub-cubic. It is then sufficient to use Chebyshev-Halley methods for schemes referencing approximate interpolant roots, without compromising the order of convergence. 
Special cases of the Chebyshev-Halley family include Chebyshev's method ($\beta=0$), Halley's method ($\beta=\frac12$) and the super-Halley method ($\beta=1$). 
%Some special cases of the Chebyshev-Halley family are presented in \cref{tbl:ChebyshevHalley}. 
The super-Halley method is generally recommended despite the leading errors for the interpolation-based schemes being $\beta$-independent, as certain next-to-leading order error terms then vanish. 
Note though that higher-order local methods~\cite{schroder1870unendlich, schroder1998infinitely} should otherwise be applied if the interpolant uses derivatives beyond first order.

%Chebyshev ($\beta=0$), Halley ($\beta=\frac12$), super-Halley ($\beta=1$)

%\begin{table}[!htbp]
%	\centering
%	%\tiny
%	%\scriptsize
%	%\footnotesize
%	\small
%	%\normalsize
%	\begin{tabular}{c@{\hspace{7mm}}c@{\hspace{7mm}}c}% p{18mm}p{14mm}p{18mm}}
%		\toprule
%		Chebyshev & Halley & super-Halley \\ \midrule
%		$\beta=0$ & $\beta=\frac12$ & $\beta=1$ \\ \bottomrule
%	\end{tabular}
%	\vspace{2mm}
%	\caption{Special cases of the Chebyshev-Halley family~\cite{gutierrez1997family}. }
%	\label{tbl:ChebyshevHalley}
%\end{table}

%\vspace{-4mm}
%\vfill\pagebreak
\section{Univariate optimisation methods}

For optimisation problems, it is first noted that interpolation-based methods derived from inverse objective functions would be unreliable. Since, if the extremum is not approached from a common direction, the multi-valued nature of the inverse function in the neighbourhood of the solution would result in a poorly constructed interpolant. Furthermore, it is non-trivial to identify the exact extremum of rational interpolants for the direct function. However, local iteration schemes about a specific point of the direct function interpolant do allow for suitable approximation of stationary points. This section presents methods based on the latter framework.

Consider the following derivative-free approximation of an objective function:
\begin{align}
\phi[x] ~\approx\,
\left.
\left( \, \sum_{i=0}^{n} \, \frac{\omega_i^{} \, \phi_i^{}}{x-x_i^{}} \, \right) \! \middle/ \! \left( \, \sum_{i=0}^{n} \, \frac{\omega_i^{}}{x-x_i^{}} \, \right) \right.
\label{eq:phi}
\end{align}

\noindent
The interpolant derivative can be deduced by performing a series expansion about a given point, which can itself be expressed by a series expansion as follows,
\begin{align}
\phi_n^\prime \hspace{2mm}\approx\hspace{2mm}
-\, \frac{1}{\omega_n^{}} \, \sum_{j \neq n}^{} \, \omega_j^{} \, \frac{ \phi_n^{}-\phi_j^{} }{x_n^{} - x_j^{} } 
\hspace{2mm}=~
- \frac{1}{\omega_n^{}} \, \sum_{k=0}^{\infty} \Bigg( \, \frac{ \phi_n^{(k+1)}}{(k+1)!} \, \sum_{j \neq n}^{} \, \omega_j^{} \,\epsilon_{jn}^{k} \, \Bigg) 
\label{eq:phi1}
\end{align}

\noindent 
where $\epsilon_{jn}^{} = (x_j^{} - x_n^{})$. In order to find general constraints for the interpolation parameters that raise the order of convergence of iteration methods, the leading errors of such schemes with respect to the true solution should be considered. However, those constraints are found to take complicated forms. This paper thus restricts its scope of discussion to identifying parameters where the derivative estimate is produced to highest order.

For \cref{eq:phi1} to reproduce the true derivative at zeroth order, it is necessary that the relation $\sum_{j\neq n}^{} \omega_j^{} = -\, \omega_n^{}$ holds, or equivalently that $\sum_i^{} \, \omega_i^{} = 0$. If series expansions for $\omega_i^{}$ trivially terminate at zeroth order, the subsequent constraints for eliminating leading errors are then equivalent to \cref{eq:generalConstraints}. As discussed previously, \cref{eq:limitSolution} offers a solution to such constraints and given the extra conditions on $\omega_i^{}$ series expansions, it follows that $\omega_i^{} \propto \prod_{j\neq i}^{} (x_i^{}-x_j^{})^{-1}$ may be applied. This choice corresponds to polynomial interpolation. 
From \cref{eq:pfe}, it then follows that the interpolant derivative equals:
\begin{align}
\phi_n^{(1)} \, -\, 
(-1)^{n} ~
\frac{ \phi_n^{(n+1)} }{(n+1)!} ~
\Bigg( \, \prod_{j=0}^{n-1}\, \epsilon_{jn}^{} \, \Bigg) 
 +\,  \text{sub-dominant terms}
\end{align}

%bracket term comes from $1/\omega_n^{}$

\noindent
The leading error when applying Newton's method with the interpolant derivatives has the following form with $m=1$:
\begin{align}
\epsilon_{n+1}^{} ~=~
%x_{\ast}^{} \, +\,  
\mathcal{E}_{m,n+1}^{} \times \Bigg( \, \epsilon_{n}^{m-1}\, \prod_{i=0}^{n-1} \, \epsilon_i^{m} \, \Bigg)
 +\,  \text{sub-dominant terms}
\label{eq:errorO}
\end{align}

\noindent
where $\mathcal{E}_{m,n+1}^{}$ is a factor composed of $k^{\text{th}}$ derivatives $\phi_\ast^{(k)}$ evaluated at the true stationary point. It is postulated that the error relation generalises as above when $m$ interpolation nodes coincide, which is indeed later confirmed for cases with $m=2$. 
In order to then identify the asymptotic order of convergence, the conventional definition that $\epsilon_{n+1}^{}  \sim \epsilon_n^\ell$ can be recursively applied to \cref{eq:errorO} to find $\ell^{n+1} = (m-1)\, \ell^{n} + m\, \sum_{i=0}^{n-1} \ell^i$. %eq~\ref{eq:lsum}. 
Equivalently, the geometric sum can be re-expressed to equate $\ell^{2} = 1  + m\, (\ell - \ell^{-n}) $. 
%Equivalently, eq~\ref{eq:lsimple} follows on re-expressing the geometric sum.
%
%\begin{eqnarray}
%	\ell^{n+1} &=& (m-1)\, \ell^{n} + m\, \sum_{i=0}^{n-1} \ell^i
%	\label{eq:lsum}
%	\\[2mm]
%	\ell^{2} &=& 1  + m\, (\ell - \ell^{-n}) 
%	\label{eq:lsimple}
%	\\[-2mm] \nonumber
%\end{eqnarray}
%
\begin{table}[!htbp]
	\centering
	%\tiny
	%\scriptsize
	%\footnotesize
	\small
	%\normalsize
	\begin{tabular}{p{22mm}p{17mm}ccccc}
		%\\[1mm]
		\toprule
		%Interpolation data & Error & 1-point & 2-point & 3-point & 4-point & 5-point \\ \midrule
		Data & Error & $n=1$ & $n=2$ & $n=3$ & $n=4$& $n\to\infty$\\ \midrule
		$\{ x_i^{} , \phi_i^{} \}$ & $\epsilon_{n}^{0} \prod_{i=0}^{n-1} \epsilon_i^{1}$ &
		1.00000 & 1.32472 & 1.46557 & 1.53416 & 1.61803 \\[1mm]
		$\{ x_i^{} , \phi_i^{}, \phi_i^\prime \}$ & $\epsilon_{n}^{1} \prod_{i=0}^{n-1} \epsilon_i^{2}$ & 
		2.00000 & 2.26953 & 2.35930 & 2.39246 & 2.41421 \\[1mm]
		$\{ x_i^{} , \phi_i^{}, \phi_i^\prime, \phi_i^{\prime\prime} \}$ & $\epsilon_{n}^{2} \prod_{i=0}^{n-1} \epsilon_i^{3}$ & 
		3.00000 & 3.22069 & 3.27902 & 3.29571 & 3.30278 \\ \bottomrule
	\end{tabular}
	\vspace{2mm}
	\caption{Convergence indices ($\ell : \epsilon_{i+1}^{} \sim \epsilon_i^\ell$) for iteration methods with a leading error given by \cref{eq:errorO}.}
	\label{tbl:rootConvergenceO} % to five decimal places
\end{table}

\noindent
\Cref{tbl:rootConvergenceO} presents numerical solutions of $\ell$. % for the cases where $m=\{1,2,3\}$. 
The order of convergence accelerates as $n$ increases, tending to $\frac12 \left( m + \sqrt{4 + m^2} \right)$. However, as for the root search methods, a long memory is not needed to raise the convergence order close to its asymptotic limit. In practice, one may again be satisfied with methods using a limited history of points.

Note that for problems where algorithmic differentiation~\cite{griewank1989automatic} is being performed to calculate gradients, the objective function is also necessarily calculated. However, many standard optimisation methods neglect the information associated with objective function values. It is inferred here that the asymptotic order convergence for methods with memory of $\{\phi_i^{}, \phi_i^\prime\}$ is $2.42$, compared to $1.62$ for the secant method. Given that the same information is calculated for both, the methods with memory converge asymptotically $1.83~(\log[2.42]/\log[1.62])$ times faster. There may also be further advantages if the convergence phase is entered faster by using non-local information.

Again, the relative efficiency of derivative-free methods and higher-derivative methods depends on the computational costs of obtaining higher derivatives. If each derivative takes a similar or longer time as $\phi[x]$ to be calculated, the derivative-free methods are generally most efficient. However, if line search optimisation is being performed as a subroutine of multi-dimensional optimisation, it will remain important to calculate gradients for directional information at certain steps. The interpolation-based approach can be adapted to consider mixed interpolation conditions, but this paper focuses on methods where the same type of information is available at each point.

\Cref{tbl:opt0} presents an optimisation scheme based on applying Newton's method to derivative-free interpolants given by \cref{eq:phi}. \Cref{tbl:opt1}  presents an optimisation scheme based on applying Chebyshev-Halley methods~\cite{gutierrez1997family} to interpolants with first derivatives, as defined below:
\begin{align}
\phi[x] ~\approx\,
\left.
\left( \, \sum_{i=0}^{n} \, \frac{\lambda_i^{} \, \phi_i^{} + (\gamma_i^{} \, \phi_i^{} + \lambda_i^{} \, \phi_i^\prime) \, (x-x_i^{}) }{(x-x_i^{})^2} \, \right) \! \middle/ \! \left( \, \sum_{i=0}^{n} \, \frac{\lambda_i^{} + \gamma_i^{} \, (x-x_i^{})}{(x-x_i^{})^2} \, \right) \right.
\end{align}

%\noindent
\vfill\pagebreak
\noindent
The leading errors in \cref{tbl:opt1} are observed to be consistent with \cref{eq:errorO}, and so the Chebyshev-Halley methods are appropriate for $n>1$ in order to not compromise the order of convergence. If Newton's method was instead applied in \cref{tbl:opt1}, the iteration scheme would be limited to quadratic convergence.

%Chebyshev ($\beta=0$), Halley ($\beta=\frac12$), super-Halley ($\beta=1$)

\vfill
%\vfill\pagebreak
%\vspace{5mm}\phantom{.}

\begin{table}[!htbp]
	\centering
	%\tiny
	%\scriptsize
	%\footnotesize
	\small
	%\normalsize
	\begin{tabular}{llp{0mm}lp{0mm}l} %p{18mm}
		\toprule
		\\[-2mm]
		\multicolumn{2}{l}{ Method: }
		&& \multicolumn{3}{l}{ $\displaystyle x_{n+1}^{} ~=~ 
			x_n^{} - \frac{\phi_n^{\prime} }{ \phi_n^{\prime\prime} }$ \hspace{2mm} : \hspace{1.5mm} 
			$
			%\begin{array}{@{}l@{}l@{}}
			\displaystyle
			\phi_n^{\prime} \hspace{2mm}\approx\, \left.
			\Bigg( \, \sum_{k \neq n}^{} \, \omega_k^{} \, \frac{ \phi_n^{}-\phi_k^{} }{x_n^{} - x_k^{} } \, \Bigg) \! \middle/ \! \Bigg( \, \sum_{k \neq n}^{} \, \omega_k^{} \, \Bigg) \right.
			%\\[4mm]
			%\displaystyle
			%\alpha_n^{-1} ~=~ 
			%\end{array}
			$ }
		\\[8mm] %\midrule %\midrule
		&&& \multicolumn{3}{l}{ $\displaystyle \phi_{n}^{\prime\prime} \hspace{2mm}\approx\, -\, 2 \left.
			\Bigg( \, \sum_{k \neq n}^{} \, \omega_k^{} \, \frac{ (\phi_n^{}-\phi_k^{}) - \phi_n^{\prime} \,(x_n^{} - x_k^{}) }{(x_n^{} - x_k^{})^2 }  \, \Bigg) \! \middle/ \! \Bigg( \, \sum_{k \neq n}^{} \, \omega_k^{} \, \Bigg) \right.
			$ }
		\\[8mm]
		%\\[-2mm]
		\multicolumn{2}{l}{ Weights: }
		&& 
		\multicolumn{3}{l}{$\displaystyle \omega_i^{} ~=~ \prod_{j\neq i}^{} \, \frac{1}{x_i^{} - x_j^{}}$
		\hspace{7mm} Error: \hspace{2mm}
		$\displaystyle \epsilon_{n+1}^{} ~\sim~ %\Bigg( 
		\frac{(-1)^{n}}{(n+1)!} ~ \frac{\phi_\ast^{(n+1)}}{ \phi_\ast^{(2)}} %\Bigg) 
		~ \Bigg( \prod_{i=0}^{n-1} \epsilon_i^{} \Bigg)$
		%\hspace{15mm} Interpolant \hspace{0mm} $\phi[x]$: (n) polynomial  
		}
		%&& %$\displaystyle \omega_i^{} ~=~ \prod_{j\neq i}^{} \, \frac{1}{f_i^{} - f_j^{}}$
		\\[5mm]
		%\multicolumn{2}{l}{ Error:} 
		%&&
		%$\displaystyle \epsilon_{n+1}^{} ~\sim~ %\Bigg( 
		%\frac{(-1)^{n}}{(n+1)!} ~ \frac{\phi_\ast^{(n+1)}}{ \phi_\ast^{(2)}} %\Bigg) 
		%~ \Bigg( \prod_{i=0}^{n-1} \epsilon_i^{} \Bigg)$
		%\\[5mm]
		\bottomrule
	\end{tabular}
	\vspace{3mm} %\small 
	\caption{Optimisation scheme based on derivative-free interpolation.}
		%(error notation defined in eq~\ref{eq:notationfi}).}
	\label{tbl:opt0}
\end{table}

%HFO with memory converges to machine precision at iterate i=5 (6th point)
%Error magnitudes for optimisation of $\phi[x] = \sin x - x^2/2$, 
%where the initial points were set to be common between methods of ...
%(using Picard iteration with $x_0^{} = 3$)

%\vfill\pagebreak
\vfill

\begin{table}[!htbp]
	\centering
	%\tiny
	%\scriptsize
	%\footnotesize
	\small
	%\normalsize
	\begin{tabular}{@{}llp{0mm}lp{1mm}l@{}} %p{18mm}
		\toprule
		\\[-2mm]
		\multicolumn{2}{l}{ Method: }
		&& \multicolumn{3}{l}{ $\displaystyle x_{n+1}^{} ~=~ 
			%x_{n}^{} - \alpha_n^{} \, \phi_{n}^{\prime} $ %\hspace{2mm} : \hspace{1mm}
			x_n^{}
			- \left( \frac
			{ \phi_n^{\prime\prime \,2} + (\frac12 - \beta ) \, \phi_n^{\prime} \, \phi_n^{\prime\prime\prime} }
			{ \phi_n^{\prime\prime \,2} - \beta \, \phi_n^{\prime} \, \phi_n^{\prime\prime\prime}  } \right) 
			\frac{ \phi_n^{\prime} }{ \phi_n^{\prime\prime} }
			%\hspace{12mm} %: \hspace{1mm} 
			\hfill
			\begin{array}{@{}l@{}}
				\beta ~\text{free parameter} \\[0.5mm]
				\beta = 1 ~\text{recommended}
			\end{array}
			$
		}
		\\[8mm] %\midrule %\midrule
		\multicolumn{6}{l}{ %\hspace{1mm}
			$\displaystyle 
			%\alpha_{n}^{-1} 
			\phi_i^{\prime\prime\phantom{\prime}} %\hspace{2mm}
			\,\approx\,
			-\, \frac{2}{\lambda_i^{}}
			\Bigg( \gamma_i^{} \, \phi_i^{\prime} 
			+\sum_{k \neq i}^{} \Bigg(
			\frac
			{  \gamma_k^{} \, (\phi_i^{} - \phi_k^{}) - \lambda_k^{}\, \phi_k^{\prime} } % (x_i^{} - x_j^{}) 
			{ (x_i^{} - x_k^{})^{} }
			+
			\frac
			{ \lambda_k^{} \, (\phi_i^{} - \phi_k^{})  }
			{ (x_i^{} - x_k^{})^2 }
			\Bigg) \!\!
			\Bigg)
			$ }
		\\[8mm]
		\multicolumn{6}{l}{ %\hspace{1mm}
		$\displaystyle 
		\phi_i^{\prime\prime\prime}
		\,\approx\,
		-\, \frac{6}{\lambda_i^{}}
		\Bigg( \frac{\gamma_i^{}\, \phi_i^{\prime\prime} }{2}
		+ \sum_{k\neq i}^{} \Bigg( 
		\frac{ \gamma_k^{} \, \phi_i^\prime  }{(x_i^{}-x_k^{})^{}}
		- \frac{ \gamma_k^{}\, (\phi_i^{} - \phi_k^{}) - \lambda_k^{}\, (\phi_i^\prime + \phi_k^\prime)   }{(x_i^{}-x_k^{})^2}
		- \frac{2\, \lambda_k^{}\, (\phi_i^{} - \phi_k^{}) }{(x_i^{}-x_k^{})^3}
		\Bigg) \!\! \Bigg)
		$  }
	\\[8mm]
		%\\[-2mm]
		\multicolumn{2}{l}{ Weights: }
		&& $\displaystyle \lambda_i^{} ~=~ \prod_{j\neq i}^{} \, \frac{1}{(x_i^{} - x_j^{})^2}$
		&& %\hspace{1mm}
		$\displaystyle \gamma_i^{} ~=~ -\, 2\, \lambda_i^{} \, \sum_{j\neq i}^{} \, \frac{1}{x_i^{} - x_j^{}}$
		\\[5mm] \midrule \midrule
		%n+1 &  
		%&& Leading error term $(\epsilon_{n+1}^{})$ 
		%&&  \\ \midrule
		\\[-2mm]
		\multicolumn{3}{l}{Error:} %Leading errors:
		& \multicolumn{3}{l}{ %\hspace{10mm}
			$\displaystyle \epsilon_{2}^{} ~\sim~ 
			-\, \frac{2\, \phi_\ast^{(4)}}{4!\, \phi_\ast^{(2)}} \, \epsilon_{0}^{2} \, \epsilon_1^{}
			\,+\,\frac{\phi_\ast^{(3)}}{2\, \phi_\ast^{(2)}} \, \epsilon_1^2 $
			\hfill
			($\beta$-independent) \hspace{5mm}
		}
		\\[7mm]
		&&& \multicolumn{3}{l}{ %\hspace{10mm}
			$\displaystyle \epsilon_{3}^{} ~\sim~ -\, \frac{2\, \phi_\ast^{(6)} }{6!\, \phi_\ast^{(2)}}\, \epsilon_0^{2}\, \epsilon_1^{2}\, \epsilon_2^{} $
		}
		\\[7mm]
		&&& \multicolumn{3}{l}{ %\hspace{10mm}
			$\displaystyle \epsilon_{4}^{} ~\sim~ -\, \frac{2\, \phi_\ast^{(8)} }{8!\, \phi_\ast^{(2)}}\, \epsilon_0^{2}\, \epsilon_1^2\, \epsilon_2^{2}\, \epsilon_3^{} $
		}
		\\[5mm] \bottomrule
	\end{tabular}
	\vspace{3mm} %\small 
	\caption{Optimisation scheme based on interpolation with first derivatives.}
		%(error notation defined in eq~\ref{eq:notationfi}).}
	\label{tbl:opt1}
\end{table}

\vfill\pagebreak
\section{Multivariate optimisation methods}

%\begin{eqnarray}
%	\phi[\mathbf{x}_i^{} + \boldsymbol{\epsilon}_i^{}]
%	&=&
%	\phi_i^{} + g_i^\alpha \, \epsilon_i^\alpha + \frac12\, H_i^{\alpha\beta} \epsilon_i^\alpha\, \epsilon_i^\beta
%	+ \mathcal{O} [\epsilon_i^3]
%\end{eqnarray}

%\begin{eqnarray}
%	\left\| \mathbf{x}-\mathbf{x}_i^{} \right\|_{M_i^{}}^2
%	%&\xrightarrow{\phantom{12}}&
%	&=&
%	\left( \mathbf{x}-\mathbf{x}_i^{} \right)^T
%	\mathbf{M}_i^{} \,
%	\left( \mathbf{x}-\mathbf{x}_i^{} \right)^{}
%\end{eqnarray}

A framework for multivariate optimisation is now proposed. However, systematic rules for parameter selection that ensure good convergence properties have not been identified. It is expected though that the univariate schemes should form a special case of suitable multivariate schemes.

The following approximation is differentiable if $\lambda_i^{}$ are non-zero parameters,
\begin{align}
\phi[\mathbf{x}]
~\approx\,
\left.
\Bigg( \, \sum_{i=0}^{n} \, \frac{\lambda_i^{} \, \phi_i^{} + ( \phi_i^{}\, \mathbf{k}_i^{}  +  \lambda_i^{} \, \mathbf{g}_i^{} )^T  (\mathbf{x}-\mathbf{x}_i^{}) }{ \left\| \mathbf{x}-\mathbf{x}_i^{} \right\|_{\mathbf{M}_i^{}}^2 } \, \Bigg) \! \middle/ \! \Bigg( \, \sum_{i=0}^{n} \, \frac{\lambda_i^{} + \mathbf{k}_i^T (\mathbf{x}-\mathbf{x}_i^{}) }{ \left\| \mathbf{x}-\mathbf{x}_i^{} \right\|_{\mathbf{M}_i^{}}^2 } \, \Bigg) \right.
\label{eq:phiNd}
\end{align}

\noindent
where $\mathbf{g}_{i}^{}$ is the gradient of the objective function at $\mathbf{x}_i^{}$, and $\{ \lambda_i^{}, \mathbf{k}_i^{} \}$ are free parameters. The notation of the norm terms is clarified below, where the $\mathbf{M}_i^{}$ metrics should be sign-definite in order for the interpolation to be well-behaved.
\begin{align}
%\left\| \mathbf{v} \right\|_{M}^2 &=& \mathbf{v}^T \mathbf{M} ~\mathbf{v}
\left\| \boldsymbol{\Delta} \right\|_{\mathbf{M}}^2 ~=~ 
\boldsymbol{\Delta}^{\!T} \mathbf{M} \, \boldsymbol{\Delta}
\end{align}

\noindent
On performing a series expansion about $\mathbf{x}_i^{}$, the interpolant derivatives are found:
%
%\vspace{0mm}
\begin{align}
%\nonumber \\[-3mm]
\left. \frac{\partial \phi^{\phantom{\alpha}}}{\partial x_{\phantom{i}}^{\alpha\phantom{\mu\!\!\!\!}}} \right|_{\mathbf{x}=\mathbf{x}_i^{}}
&~=~
g_i^{\alpha\phantom{\mu\!\!\!\!}}
\\[2mm]
H_i^{\mu\nu}
~=~
\left. \frac{\partial^2 \phi }{\partial x_{\phantom{i}}^\mu \, \partial x_{\phantom{i}}^{\nu\phantom{\mu\!\!\!\!}} } \right|_{\mathbf{x}=\mathbf{x}_i^{}}
%H_{i}^{\mu\nu}
%&\approx&
&~=~
- \frac{ s_i^{} \,  M_i^{\mu \nu}  + g_i^{\mu} \, k_i^{\nu\phantom{\mu\!\!\!\!}} + k_i^{\mu} \, g_i^{\nu\phantom{\mu\!\!\!\!}}  }{\lambda_i^{}}
\label{eq:Hessian}
\end{align}
%
%\vspace{5mm}
\begin{align}
\text{where} \hspace{2mm}
s_i^{}  ~=~
2\, \sum_{j \neq i}^{} 
\Bigg(
\frac{ \lambda_j^{}\,  (\phi_i^{} - \phi_j^{} ) + ( (\phi_i^{} - \phi_j^{}) \, \mathbf{k}_{j}^{} - \lambda_j^{} \, \mathbf{g}_j^{} )^T (\mathbf{x}_i^{} -\mathbf{x}_j^{})  }{ 
	\left\| \mathbf{x}_i^{} -\mathbf{x}_j^{} \right\|_{\mathbf{M}_j^{}}^2
	%\left( \mathbf{x}_i^{} -\mathbf{x}_j^{} \right)^{\!T}
	%\mathbf{M}_j^{} 
	%\left( \mathbf{x}_i^{} -\mathbf{x}_j^{} \right)^{} 
	}
\Bigg)
%\\[-2mm] \nonumber 
\end{align}

%\vspace{1mm}
\noindent
The interpolant Hessian $\mathbf{H}_i^{}$ possesses a sufficient number of degrees of freedom so that it could tend to the true Hessian. Depending on the dimension, and number of interpolation points, it may be convenient to decompose $\mathbf{M}_i^{}$ into low rank terms plus a term proportional to the identity matrix. However, it is not obvious what basis to choose for such decompositions.

The inverse Hessian can be calculated by repeated use of the Sherman-Morrison formula, resulting in the form  $\mathbf{H}_i^{-1} = \mathbf{M}_{i}^{-1}\, \mathbf{W}_i^{}~ \mathbf{M}_{i}^{-1} $ with the following definitions:
\begin{align}
%\nonumber \\[-3mm]
%\left( M_i^{} \, H_i^{-1} M_i^{}  \right)^{\mu\nu} \!\!
W_i^{\mu\nu} 
%&\approx& \!\!
\,=\,
%\!\! &=& \!\!
-\, \frac{\lambda_i^{}}{s_i^{}}
\left( 
M_i^{\mu\nu}
%\delta^{\mu\nu}
%\left( M_i^{-1} \right)^{\mu\nu}
+\frac
{ g_i^{\mu} g_i^{\nu\phantom{\mu\!\!\!\!}} \left\| \mathbf{k}_i^{} \right\|_{\mathbf{M}_i^{-1}}^2 
	+ k_i^{\mu} k_i^{\nu\phantom{\mu\!\!\!\!}} \left\| \mathbf{g}_i^{} \right\|_{\mathbf{M}_i^{-1}}^2
	-( g_i^{\mu} k_i^{\nu\phantom{\mu\!\!\!\!}} + k_i^{\mu} g_i^{\nu\phantom{\mu\!\!\!\!}} ) \,t_i^{} }
{ t_i^2 \,-\, \left\| \mathbf{k}_i^{} \right\|_{\mathbf{M}_i^{-1}}^2 \left\| \mathbf{g}_i^{} \right\|_{\mathbf{M}_i^{-1}}^2 }
\right) %( s_i^{} + k_i^{} .g_i^{} )^2
%~\phantom{.}
\end{align}
\begin{align}
	%\mathbf{H}_i^{-1} 
	%~\approx~ \mathbf{M}_{i}^{-1}\, \mathbf{W}_i^{}~ \mathbf{M}_{i}^{-1} 
	%\hspace{10mm}
	\text{where}\hspace{2mm} t_i^{} ~=~ \mathbf{g}_i^T \, \mathbf{M}_i^{-1} \mathbf{k}_i^{} + s_i^{} 
	%\\[-3mm] \nonumber 
\end{align}

\noindent
On applying Newton's method to the interpolant in \cref{eq:phiNd}, the following iteration scheme is then established:
\begin{align}
%\boxed{\hspace{1mm}
%	\begin{array}{c}
%	\\[-3mm]
%	\displaystyle
	\mathbf{x}_{n+1}^{}
	~=~
	\mathbf{x}_n^{} - \mathbf{H}_n^{-1} \, \mathbf{g}_n^{} 
	%\\[4mm]
	~=~
	\mathbf{x}_n^{} + \lambda_n^{}\, \mathbf{M}_n^{-1} \left( 
	\frac
	{ t_n^{}\, \mathbf{g}_n^{} - \mathbf{k}_n^{} \left\| \mathbf{g}_n^{} \right\|_{\mathbf{M}_n^{-1}}^2 }
	{ t_n^2 - \left\| \mathbf{k}_n^{} \right\|_{\mathbf{M}_n^{-1}}^2 \left\| \mathbf{g}_n^{} \right\|_{\mathbf{M}_n^{-1}}^2 } 
	\right)
%	\\[4mm]
%	\end{array}
%}
\end{align}

\noindent
Although the above formula is motivated by interpolation, it may be helpful to simply assume the Hessian form in \cref{eq:Hessian} and then require that certain higher derivatives are constant. Constraints can then be set for a subset of parameters. However, the question of how to ensure good convergence properties still remains.

\vfill\pagebreak
\section{Multivariate root search methods}

The approaches considered for univariate root search can also be applied in the multivariate case. 
However, as for the proposed multivariate optimisation framework, systematic rules for parameter selection that ensure good convergence properties have not been identified.

%approaches can be based on either exact or approximate roots of interpolants for the inverse or direct function, similarly to the univariate frameworks.

To derive methods based on inverse interpolation, consider the following relation:
\begin{align}
%\small 
\mathbf{x} [\mathbf{f}]
~ \approx\,
\Bigg( \, \sum_{i=0}^{n}  \frac{ \boldsymbol{\Lambda}_i^{} +  \boldsymbol{\Gamma}_i^{}  \, (\mathbf{f} - \mathbf{f}_i^{} ) }{\left\| \mathbf{f} - \mathbf{f}_i^{} \right\|_{\mathbf{M}_i^{}}^2} \, \Bigg)^{\!\!\! -1} 
\Bigg( \, \sum_{i=0}^{n}  \frac{ \boldsymbol{\Lambda}_i^{} \, \mathbf{x}_i^{} + ( \boldsymbol{\Gamma}_i^{} \, \mathbf{x}_i^{} + \boldsymbol{\Lambda}_i^{} \, \mathbf{J}_{i}^{-1})  (\mathbf{f} - \mathbf{f}_i^{} ) }{\left\| \mathbf{f} - \mathbf{f}_i^{} \right\|_{\mathbf{M}_i^{}}^2} \, \Bigg)
%\phantom{.}
\label{eq:xfmulti}
\end{align}

\noindent
where $\mathbf{J}_i^{}$ are (possibly approximated) Jacobian matrices, $\boldsymbol{\Lambda}_i^{}$ are free matrices, and $\boldsymbol{\Gamma}_i^{}$ free tensors 
such that %$\boldsymbol{\Gamma}_i^{}  \mathbf{x}_i^{}  \mathbf{f}_i^{}  = \Gamma_i^{\mu \nu \rho}  x_i^{\nu\phantom{\mu\!\!\!\!}}  f_i^{\rho}$ and 
$\boldsymbol{\Gamma}_i^{} \, \mathbf{f}_i^{}  = \Gamma_i^{\mu \nu \rho}  f_i^{\rho}$ and
$\boldsymbol{\Gamma}_i^{}  \mathbf{x}_i^{}   \mathbf{f}_i^{}  = \Gamma_i^{\mu \nu \rho}  x_i^{\nu \phantom{\mu\!\!\!\!}} f_i^{\rho}$. 
Note that for the interpolant to be differentiable, the $\boldsymbol{\Lambda}_i^{}$ matrices must be non-zero. 
On notation conventions, it is assumed in this section that pairs of (Greek) dimension indices are summed over when appearing in a given term. 
The following iteration formula is then found on setting $\mathbf{f}=\mathbf{0}$:
\begin{align}
%\boxed{\hspace{1mm}
%	\begin{array}{c}
%	\\[-3mm]
	\mathbf{x}_{n+1}^{}
	~=\,
	\Bigg( \,
	\sum_{i=0}^{n} \,
	\frac
	{ \boldsymbol{\Lambda}_i^{} - \boldsymbol{\Gamma}_i^{} \, \mathbf{f}_i^{} }
	{\mathbf{f}_i^{T\phantom{-}\!\!\!\!} \mathbf{M}_i^{} \, \mathbf{f}_i^{\phantom{-}\!\!\!\!}}
	\, \Bigg)^{\!\!\! -1}
	\Bigg( \,
	\sum_{i=0}^{n} \,
	\frac
	{ \boldsymbol{\Lambda}_i^{} \left( \mathbf{x}_i^{} - \mathbf{J}_i^{-1} \mathbf{f}_i^{\phantom{-}\!\!\!\!} \right) - \boldsymbol{\Gamma}_i^{} \, \mathbf{x}_i^{} \, \mathbf{f}_i^{} }
	{\mathbf{f}_i^{T\phantom{-}\!\!\!\!} \mathbf{M}_i^{} \, \mathbf{f}_i^{\phantom{-}\!\!\!\!}}
	\, \Bigg)
%	\\[4mm]
%	\end{array}
%}
\label{eq:xmulti}
\end{align}

\noindent
The above iteration formula can also be derived by approximating the direct function in \cref{eq:fxmulti}, given the relation for the $\mathbf{R}$ matrix function in \cref{eq:rApprox}:
\begin{align} %\mathbf{h}
%\nonumber \\[-3mm]
\mathbf{f} [\mathbf{x} ]
~=~
( \mathbf{R}_{}^{-1} [\mathbf{x}]  )
\left( \mathbf{A}\, \mathbf{x} + \mathbf{b}  \right) %\boldsymbol{\xi}
\label{eq:fxmulti}
\end{align}
\begin{align}
%\small 
( \mathbf{R} \, \mathbf{f} ) [\mathbf{x}]
~ \approx\,
\left( \, \sum_{i=0}^{n}  \frac{ \boldsymbol{\Lambda}_i^{} 
	%+  \boldsymbol{\Gamma}_i^{}   (\mathbf{x} - \mathbf{x}_i^{} ) 
	}{\left\| \mathbf{x} - \mathbf{x}_i^{} \right\|^2} 
	\right)^{\!\!\! -1} 
\left( \, \sum_{i=0}^{n}  \frac{ \boldsymbol{\Lambda}_i^{} \, ( \mathbf{R}_i^{} \, \mathbf{f}_i^{} 
	%+ ( \boldsymbol{\Gamma}_i^{} \, \mathbf{x}_i^{} 
	+ (\mathbf{R}_i^{\prime} \, \mathbf{f}_i^{} + \mathbf{R}_i^{} \, \mathbf{J}_i^{}) \, (\mathbf{x} - \mathbf{x}_i^{} )  )
 }{\left\| \mathbf{x} - \mathbf{x}_i^{} \right\|^2}  \right)
\label{eq:rApprox}
%\\[-2mm] \nonumber 
\end{align}

\noindent
where $\mathbf{R}_i^{\prime}$ is the derivative of the $\mathbf{R}$ matrix function at $\mathbf{x}_i^{}$, and $(\mathbf{R}_i^\prime \, \mathbf{f}_i^{} )_{\phantom{i}}^{\mu\nu} = R_i^{\prime \mu\sigma\nu} f_i^{\sigma\phantom{\mu}\!\!\!\!}$. 
If the conditions below hold, \cref{eq:fxmulti} and \cref{eq:rApprox} are guaranteed to be consistent.
\begin{align}
\begin{array}{rcl@{\hspace{6mm}}l}
\mathbf{A}
\! &=& \!
\mathbf{R}_i^{\prime} \, \mathbf{f}_i^{} + \mathbf{R}_i^{} \, \mathbf{J}_i^{}
& \forall \, i
\\[3mm]
\mathbf{b}
\! &=& \!
\mathbf{R}_i^{} \, \mathbf{f}_i^{}
-  ( \mathbf{R}_i^{\prime} \, \mathbf{f}_i^{} + \mathbf{R}_i^{} \, \mathbf{J}_i^{} ) \, \mathbf{x}_i^{}
& \forall \, i
\end{array}
\label{eq:rConstraint}
\end{align}

\noindent
Solutions for \cref{eq:rConstraint} are given by:
\begin{align}
%\nonumber \\[-3mm]
&
R_i^{\mu\nu} ~~~\phantom{.} \hspace{1mm}
~=~ \hspace{0mm}
%\Bigg( 
\Bigg( \!\! \Bigg( \,
\sum_{k = 0}^{n} \, 
\frac{ \Lambda_{k}^{\mu\alpha} - \Gamma_k^{\mu\alpha\gamma } f_k^\gamma  }{a_k^{}} 
%-\frac{ \Lambda_{k}^{\mu\gamma}  ( z_k^{\gamma} - z_i^{\gamma} ) \, f_i^{\sigma\phantom{\mu}\!\!\!\!} \, \tilde{M}_i^{\sigma \alpha \phantom{\mu}\!\!\!\!} }{a_k^{} \, a_i^{}} 
\, \Bigg)
- R_i^{\prime \mu\sigma\alpha} f_i^{\sigma \phantom{\mu}\!\!\!\! }
\, \Bigg) \,
%\Bigg)
\Big( J_i^{-1} \Big)^{\! \alpha\nu}
~\phantom{.}
\\[2mm]
&
R_i^{\prime \mu\sigma\nu} \hspace{3mm}
~=~ \hspace{0mm}
\Bigg( \,
\sum_{k = 0}^{n} \, 
%\frac{ \Lambda_{k}^{\mu\gamma}  ( z_k^{\gamma} - z_i^{\gamma} ) \, \tilde{M}_i^{\sigma \nu \phantom{\mu}\!\!\!\!} }{a_k^{} \, a_i^{}} 
\frac{ \Lambda_{k}^{\mu\gamma}  ( z_k^{\gamma} - z_i^{\gamma} ) - \Gamma_k^{\mu\gamma\rho} (x_k^{\gamma}  - z_i^{\gamma} )  f_k^{\rho} }{a_k^{}} 
\, \Bigg) \,
\frac{ \tilde{M}_i^{\sigma \nu \phantom{\mu}\!\!\!\!} }{ a_i^{} } 
\\[4mm]
&
\text{where}\hspace{5mm}
a_i^{} ~=~ \mathbf{f}_i^{T\phantom{-}\!\!\!\!} \tilde{\mathbf{M}}_i^{} \, \mathbf{J}_i^{-1} \mathbf{f}_i^{\phantom{-}\!\!\!\!}
\hspace{10mm}
\mathbf{z}_i^{} ~=~ 
\mathbf{x}_i^{} - \mathbf{J}_i^{-1} \mathbf{f}_i^{\phantom{-}\!\!\!\!}
%\nonumber
%\\[-3mm] \nonumber 
\end{align}

\noindent
The $\mathbf{A}, \mathbf{b}$ factors are then equal to: 
\begin{align}
%\bar{\mathbf{K}}
\mathbf{A}
~=\,
%\phantom{- \,}
\Bigg( \,
\sum_{k=0}^{n} \,
\frac
{ \boldsymbol{\Lambda}_k^{} - \boldsymbol{\Gamma}_k^{} \, \mathbf{f}_k^{} }
{\mathbf{f}_k^{T\phantom{-}\!\!\!\!} \, \tilde{\mathbf{M}}_k^{} \, \mathbf{J}_k^{-1} \, \mathbf{f}_k^{\phantom{-}\!\!\!\!}}
\, \Bigg)
\hspace{5mm}
\mathbf{b}
~=\,
- \,
\Bigg( \,
\sum_{k=0}^{n} \,
\frac
{ \boldsymbol{\Lambda}_k^{} \left( \mathbf{x}_k^{} - \mathbf{J}_k^{-1} \mathbf{f}_k^{\phantom{-}\!\!\!\!} \right) - \boldsymbol{\Gamma}_k^{} \, \mathbf{x}_k^{} \, \mathbf{f}_k^{} }
{\mathbf{f}_k^{T\phantom{-}\!\!\!\!} \, \tilde{\mathbf{M}}_k^{} \, \mathbf{J}_k^{-1} \, \mathbf{f}_k^{\phantom{-}\!\!\!\!}}
\, \Bigg)
\end{align}

\noindent
The iteration formula in \cref{eq:xmulti} is thus reproduced by solving for $\mathbf{f} = \mathbf{0}$ given~\cref{eq:fxmulti} $(\mathbf{x}_{n+1}^{} = - \mathbf{A}^{-1} \mathbf{b})$ and defining $\mathbf{M}_i^{} = \tilde{\mathbf{M}}_i^{} \, \mathbf{J}_i^{-1}$.

Alternatively, root search methods can be constructed by applying local iteration schemes with interpolant derivatives, either of $\mathbf{x}[\mathbf{f}]$ as in \cref{eq:xfmulti}, or its counterpart form for $\mathbf{f}[\mathbf{x}]$. On convergence behaviours of the proposed approaches, naive extensions of the univariate parameter choices do not achieve convergence acceleration, and also tend to delay the convergence phase. Of course the univariate method parameters were required to obey specific constraints. A full set of corresponding constraints for multivariate schemes has not yet been determined.

% Chebyshev-Halley family~\cite{gutierrez1997family} 

%\vfill\pagebreak
\section{Conclusions}

The convergence order of univariate root search and optimisation schemes can be accelerated by re-using accumulated function information, given certain constraints identified in this paper. However, a long memory is not required in order to approach the asymptotic limits. For univariate root search, the presented derivative-free methods approach quadratic convergence and the first-derivative methods approach cubic convergence. For univariate optimisation, the derivative-free methods approach a convergence order of $1.62$ and the first-derivative methods approach an order of $2.42$. There are general performance advantages with respect to low-memory iteration methods, most notably in the case where optimisation routines calculate the objective function and gradient at each step: the full-memory methods converge asymptotically $1.8$~times faster than the secant method. For problems where the time required to calculate derivatives is similar or longer than that for function evaluations, it is stressed that the derivative-free methods are most efficient.

Frameworks to extend the iteration schemes to multivariate problems have also been proposed, but without identification of practically useful parameter choices. Further study is therefore required to answer if/how the approaches can be suitably extended to multivariate problems. There are also various other ways that this work can be extended. Although not listed in this paper, iteration formulae with mixed interpolation conditions can be defined. Furthermore, problems with non-simple roots or stationary points with null second derivatives require adaptations to the interpolation-based schemes presented. Another restriction within this paper was that only `greedy' approaches were considered, where the leading error of the next step was required to be maximally suppressed. It may be expected though that the error tolerance will not be met in the subsequent step (given knowledge of the convergence rates), and so multi-point methods with memory should then be favoured. 

For final emphasis, the univariate iteration schemes presented in this paper are advocated for common use. The associated parameters are defined analytically, and so require no intermediate calibration operations. The formulae can also be applied to an arbitrary history of points. Primarily though, performance advantages are achieved when the iteration evaluations form a negligible component of the computation efforts.

%Some care must be taken to manage numerical precision errors, but the implementation can be performed efficiently. 

%numerical issues?
%analytic acceleration parameters, easy to update
%stochastic applications, with uncertainty estimation
%robust to have non-local information, attraction basins?
%invariant under affine co-ordinate transformations

%tuning the iteration efficiency for specific problems
%transform function, co-ordinate spaces

%other interpolation/approximation approaches
%approximate the gradient, or Jacobian

%\vfill\pagebreak
\bibliographystyle{unsrt}

\end{document}